\documentclass[conference,compsoc]{IEEEtran}

\newtheorem{theorem}{\indent Theorem}
\newtheorem{remark}{\indent Remark}

\usepackage[intlimits]{amsmath}
\usepackage{amssymb}

\ifCLASSOPTIONcompsoc
  \usepackage[nocompress]{cite}
\else
  \usepackage{cite}
\fi

\ifCLASSINFOpdf
  \usepackage[pdftex]{graphicx}
\else
\fi

\usepackage{amsmath}

\begin{document}
\def \op {\overline{p}}
\def \oP {\overline{P}}
\def \oz {\overline{\zeta}}
\def \ou {\overline{u}}
\def \oM {\overline{M}}
\def \oG {\overline{G}}
\def \tM {\tilde{M}}
\def \eps {\varepsilon}
\def \V {\mathrm{Var}}
\def \mE {\mathbb{E}}
\def \heta {\hat{\eta}}
\def \hzeta {\hat{\zeta}}
\def \hxi {\hat{\xi}}
\def \hp {\hat{p}}
\def \hL {\hat{\mathrm{L}}}
\def \L {\mathrm{L}}
\def \hD {\hat{D}}
\def \d {\delta}
\def \mP {\mathcal{P}}
%
\title{Two-Armed Bandit Problem, Data Processing,\\ and Parallel Version of the Mirror Descent Algorithm}

\author{\IEEEauthorblockN{Alexander Kolnogorov}
\IEEEauthorblockA{Yaroslav-the-Wise\\ Novgorod State University\\
B.St-Petersburgskaya Str, 41\\
Velikiy Novgorod, Russia, 173003\\
Email: Alexander.Kolnogorov@novsu.ru} \and
\IEEEauthorblockN{Alexander Nazin}
\IEEEauthorblockA{V.A.Trapeznikov Institute of Control Sciences \\
Russian Academy of Sciences\\
65, Profsoyuznaya str.\\
117997 Moscow, Russia\\
Email: nazine@ipu.ru}
\and \IEEEauthorblockN{Dmitry Shiyan}
\IEEEauthorblockA{Yaroslav-the-Wise\\ Novgorod State University\\
B.St-Petersburgskaya Str, 41\\
Velikiy Novgorod, Russia, 173003\\
Email: dsqq-tm@ya.ru}}

\maketitle

\begin{abstract}
We consider the minimax setup for the two-armed bandit problem as
applied to data processing if there are two alternative processing
methods available with different a priori unknown efficiencies.
One should determine the most effective method and provide its
predominant application. To this end we use the mirror descent
algorithm (MDA). It is well-known that corresponding minimax risk
has the order $N^{1/2}$ with $N$ being the number of processed
data. We improve significantly the theoretical estimate of the
factor using Monte-Carlo simulations. Then we propose a parallel
version of the MDA which allows processing of data by packets in a
number of stages. The usage of parallel version of the MDA ensures
that total time of data processing depends mostly on the number of
packets but not on the total number of data. It is quite
unexpectedly that the parallel version behaves unlike the ordinary
one even if the number of packets is large. Moreover, the parallel
version considerably improves control performance because it
provides significantly smaller value of the minimax risk. We
explain this result by considering another parallel modification
of the MDA which behavior is close to behavior of the ordinary
version. Our estimates are based on invariant descriptions of the
algorithms. All estimates are obtained by Monte-Carlo simulations.
\par
It's worth noting that parallel version performs well only for
methods with close efficiencies. If efficiencies differ
significantly then one should use the combined algorithm which at
initial sufficiently short control horizon uses ordinary version
and then switches to the parallel version of the MDA.
\end{abstract}

\noindent {\small \bf \emph{Keywords---} two-armed bandit problem,
control in a random environment, minimax approach, robust control,
mirror descent algorithm, parallel processing.}

\IEEEpeerreviewmaketitle

\section{Introduction}
We consider the two-armed bandit problem (see, e.g. \cite{BF},
\cite{PS}) which is also well-known as the problem of expedient
behavior in a random environment (see, e.g. \cite{Tsetlin},
\cite{Varsh}) and the problem of adaptive control in a random
environment (see, e.g. \cite{Sragovich}, \cite{NP}) in the
following setting. Let $\xi_n$, $n=1,\dots,N$ be a controlled
random process which values are interpreted as incomes, depend
only on currently chosen actions $y_n$ ($y_n \in \{1,2\}$) and
have probability distributions
\begin{gather*}
\Pr(\xi_n =1|y_n=\ell)=p_\ell,\quad \Pr(\xi_n =0|y_n=\ell)=q_\ell
\end{gather*}
where $p_\ell+q_\ell=1$, $\ell=1,2$. So, this is the so-called
Bernoulli two-armed bandit. Such bandit is described by a
parameter $\theta=(p_1,p_2)$ with the set of values
$\Theta=\{\theta:0\leq p_\ell\leq 1; \ell=1,2\}$. It is well-known
that mathematical expectation and variance of one-step income are
equal to
\begin{gather*}
m_\ell=\mE (\xi_n|y_n=\ell)=p_\ell,\
D_\ell=\V(\xi_n|y_n=\ell)=p_\ell q_\ell
\end{gather*}
\par
The goal is to maximize (in some sense) the total expected income.
Control strategy $\sigma$ at the point of time $n$ assigns a
random choice of the action $y_n$ depending on the current history
of the process, i.e. responses $x^{n-1}=x_1,\dots,x_{n-1}$ to
applied actions $y^{n-1}=y_1,\dots,y_{n-1}$:
$$
\sigma_{\ell}(y^{n-1},x^{n-1})=\Pr(y_n=\ell|y^{n-1},x^{n-1}) ,
$$
$\ell=1,2$. The most general set of strategies is denoted here by
$\Sigma_0$. However, there may be some additional restrictions on
the set of strategies. For example, one can consider the set of
strategies $\Sigma_1$ described by the mirror descent algorithm
presented below. More restrictive is additional requirement to
strategy to allow parallel processing, this is the set $\Sigma_2$
below. In the sequel we define some additional sets of strategies.
\par
If parameter $\theta$ is known then the optimal strategy should
always apply the action corresponding to the largest value of
$m_1$, $m_2$. The total expected income would thus be equal to $N
(m_1 \vee m_2)$. If parameter is unknown then the regret
\begin{equation*}
L_N(\sigma,\theta)=N (m_1 \vee m_2)-
\mE_{\sigma,\theta}\left(\sum_{n=1}^N  \xi_n\right)
\end{equation*}
describes expected losses of total income with respect to its
maximal possible value due to incomplete information. Here
$\mathbb{E}_{\sigma,\theta}$ denotes the mathematical expectation
calculated with respect to measure generated by strategy $\sigma$
and parameter $\theta$.
\par
According to the minimax approach the maximal value of the regret
on the set of parameters $\Theta$ should be minimized on the set
of strategies. The value
\begin{equation}\label{a1}
R^{(0)}_N(\Theta)=\inf_{\Sigma_0} \sup_{\Theta} L_N(\sigma,\theta)
\end{equation}
is called the minimax risk corresponding to the most general set
of strategies $\Sigma_0$ and the optimal strategy $\sigma_0^M$ is
called the minimax strategy. Application of the minimax strategy
ensures that the following inequality holds
\begin{equation*}
L_N(\sigma_0^M,\theta) \le R^{(0)}_N(\Theta)
\end{equation*}
for all $\theta \in \Theta$ and this implies robustness of the
control.
\par
The minimax approach to the problem was proposed by H.~Robbins
in~\cite{Robbins}. This article caused a significant interest to
considered problem. It was shown in \cite{FZ} that explicit
determination of the minimax strategy and minimax risk is
virtually impossible already for $N>4$. However, the following
asymptotic minimax theorem was proved by W.~Vogel in \cite{Vogel}.
\begin{theorem} The following estimates of the minimax risk~\eqref{a1}
hold as $N \to \infty$ for
Bernoulli two-armed bandit:
\begin{equation}
0.612  \leq (DN)^{-1/2}R^{(0)}_N(\Theta) \leq 0.752 \label{a2}
\end{equation}
with $D=0.25$ being the maximal possible variance of one-step
income.
\end{theorem}
Presented here the lower estimate was obtained in \cite{Bather}.
The upper estimate was obtained in~\cite{Vogel} for the following
strategy.
\par
{\bf Thresholding strategy.} \emph{ Use actions turn-by-turn until
the absolute difference between total incomes for their
applications exceeds the value of the threshold $\alpha
(DN)^{1/2}$ or the control horizon expires. If the threshold has
been achieved and the control horizon has not expired then at the
residual control horizon use only the action corresponding to the
larger value of total initial income.}
\par
The optimal value of $\alpha$ is $\alpha\approx 0.584$ and the
maximal value of the regret corresponds to $|p_1-p_2|\approx 3.78
(D/N)^{1/2}$ with additional requirement that $p_1$, $p_2$ are
close to 0.5.
\par
This approach is generalized in \cite{Koln1, Koln2} for Gaussian
(or Normal) two-armed bandit, i.e. described by the probability
distribution density of incomes
$$
f_D(x|m_\ell)=(2\pi D)^{-1/2} \exp
\left\{-(x-m_\ell)^2/(2D)\right\}
$$
if $y_n=\ell$ ($\ell=1,2$). Gaussian two-armed bandit can be
described by a vector parameter $\theta=(m_1, m_2)$. The set of
parameters is assumed to be the following
$$\Theta=\{\theta:|m_1-m_2| \leq 2C\},$$
where $0<C< \infty$. Restriction $C < \infty$ ensures the
boundedness of the regret on $\Theta$.
\par
In \cite{Koln1, Koln2}, according to the main theorem of the
theory of games the minimax risk for Gaussian two-armed bandit is
sought for as Bayesian one corresponding to the worst-case prior
distribution for which Bayesian risk attains its maximal value.
The Bayesian approach allows to write recursive Bellman-type
equation for numerical determination of both Bayesian strategy and
Bayesian risk. However, a direct application of the main theorem
of the theory of games is virtually impossible because of its high
computational complexity. Therefore, at first a description of the
worst-case prior distribution is done. It is shown that the
worst-case prior is symmetric and asymptotically uniform and this
allows significantly to symplify the Bellman-type equation.
In~\cite{Koln2} the estimates \eqref{a2} are improved as follows.
\par
\begin{theorem} The following estimate of the minimax risk~\eqref{a1} holds for
Gaussian two-armed bandit
\begin{equation}
\lim_{N \to \infty}(DN)^{-1/2}R^{(0)}_N(\Theta)=r_0 \label{a3}
\end{equation}
with $r_0 \approx 0.637$.
\end{theorem}
\begin{remark}
In \cite{Koln1, Koln2} the case $D=1$ is considered. However, all
reasonings can be easily extended to distributions with arbitrary
$D$.
\end{remark}
Let's explain the choice of Gaussian distribution of incomes. We
consider the problem as applied to group control of processing a
large amount of data. Let $N=TM$ data items be given that can be
processed using either of the two alternative methods. The result
of processing of the $n$-th item of data is $\xi_n = 1$ if
processing is successful and $\xi_n = 0$ if it is unsuccessful.
Probabilities $\Pr(\xi_n=1|y_t=\ell)=p_\ell$, $\ell=1,2$ depend
only on selected methods (actions). Let's assume that $p_1$, $p_2$
are close to $p$ ($0<p<1$). We partition the data into $T$ packets
of $M$ data in each packet and use the same method for data
processing in the same packet. For the control, we use the values
of the process $\eta_t=M^{-1/2} \sum_{n=(t-1)M+1}^{tM} \xi_n$,\
$t=1,\dots,T$. According to the central limit theorem probability
distributions of $\eta_t$, $t=1,\dots,T$ are close to Gaussian and
their variances are close to $D=p(1-p)$ as in considered setting.
\par
Note that data in the same packet may be processed in parallel. In
this case, the total time of data processing depends on the number
of packets rather than on the total number of data. However, there
is a question of losses in the control performance as the result
of such aggregation.  It was shown in~\cite{Koln1, Koln2} that if
$T$ is large enough (e.g. $T \ge 30$) then parallel control is
close to optimal. Therefore, say $30000$ items of data can be
processed in 30 stages by packets of 1000 data with almost the
same maximal losses as if the data were processed optimally
one-by-one. However, one should ensure the closeness of $p_1$,
$p_2$ in this case. Otherwise parallel processing causes large
losses at the initial stage of the control when both actions are
applied turn-by-turn. In~\cite{Koln1, Koln2} this requirement is
discussed in more details and an adaptive algorithm which is
optimal  for both close and distant $p_1$, $p_2$ is proposed.
\par
\begin{remark}
The estimates \eqref{a2} can be easily extended to Gaussian
two-armed bandit with a glance that the maximal value of the
regret corresponds to $|m_1-m_2|\approx 3.78 (D/N)^{1/2}$ in this
case. In particular, this implies that thresholding strategy
allows parallel processing. The estimate \eqref{a3} can be in turn
extended to Bernoulli two-armed bandit by usage of parallel
processing of data.
\end{remark}
\par
\begin{remark}
Parallel control in the two-armed bandit problem was first
proposed for treating a large group of patients by either of the
two alternative drugs with different unknown efficiencies.
Clearly, the doctor cannot treat the patients sequentially
one-by-one. Say, if the result of the treatment will be manifest
in a week and there is a thousand of patients, then one-by-one
treatment would take almost twenty years. Therefore, it was
proposed to give both drugs to sufficiently large test groups of
patients and then the more effective one to give to the rest of
them. As the result, the entire treatment takes two weeks. Note
that the two-armed bandit problem, as applied to medical trials,
was usually considered in Bayesian setting and for sufficiently
small number of stages (two, and sometimes three or four treatment
stages). So, the results depend on the prior which is often
specifically chosen and the control quality is less than for
sufficiently large number of stages. The discussion and
bibliography of the problem can be found, for example,
in~\cite{Lai}.
\end{remark}
\par
There are some different approaches to robust control in the
two-armed bandit problem, see, e.g. \cite{NP, Lugosi, UN, GNS}. In
these articles, stochastic approximation method and mirror descent
algorithm are used for the control. Instead of minimax risk, the
authors often consider the equivalent attitude called the
guaranteed rate of convergency. The order of the minimax risk for
these algorithms is $N^{1/2}$ or close to $N^{1/2}$. However, more
precise estimates were not presented. The versions for parallel
processing were not proposed as well.
\par
The goal of this paper is to investigate the mirror descent
algorithm  (MDA) for the two-armed bandit problem proposed
in~\cite{UN}. For this algorithm the minimax risk has the order
$N^{1/2}$ and theoretical estimate of the factor (or normalized
minimax risk) is $r_1 \le 4.710$. We improve this estimate by
Monte-Carlo simulations as $r_1 \le 2.0$. Then we propose a
parallel version of the algorithm which partitions application of
actions in the packet in proportion to corresponding
probabilities. For this parallel version of the MDA, we obtain
invariant description which does not depend on the size of the
packet. We show that corresponding minimax risk has the order
$N^{1/2}$ and estimate the value of the factor as $r_2 \approx
1.1$ using Monte-Carlo simulations. It is quite unexpectedly that
parallel version behaves unlike the ordinary one even if the
number of packets is sufficiently large. Moreover, it provides
significantly smaller value of the minimax risk. We explain this
result by considering another parallel version of the MDA which
partitions actions in the packet sequentially with probabilities
determined at the beginning of packet processing. This version of
the MDA behaves like the ordinary one if the number of packets is
large enough. For this version of MDA, we obtain invariant
description as well.
\par
It is important that parallel versions of the MDA perform well
only for close values of probabilities $p_1$, $p_2$. For distant
probabilities there may be significant expected losses caused by
processing of the first packet. To avoid this,  combined versions
of the MDA are proposed.  These algorithms at initial sufficiently
short stage apply the ordinary MDA and then switch to the parallel
version. These algorithms  perform well for all probabilities
$p_1$, $p_2$.
\par
The structure of the paper is the following. In Section~\ref{S2}
we present the description of the algorithm from~\cite{UN} and
improve the estimate of the minimax risk by Monte-Carlo
simulations. In Section~\ref{S3} we propose the version of this
algorithm which allows parallel processing and propose the
invariant description of the algorithm. In Section~\ref{S4} we
propose another parallel version of the MDA which behaves like the
ordinary algorithm. Combined algorithms are presented in
Section~\ref{S5}. Section~\ref{Con} contains a conclusion. Note
that some results were presented in~\cite{Koln3}.

\section{Description and Properties of the MDA for Bernoulli Two-Armed Bandit } \label{S2}

\begin{figure}[h]
\centering
\includegraphics[width=3.2in]{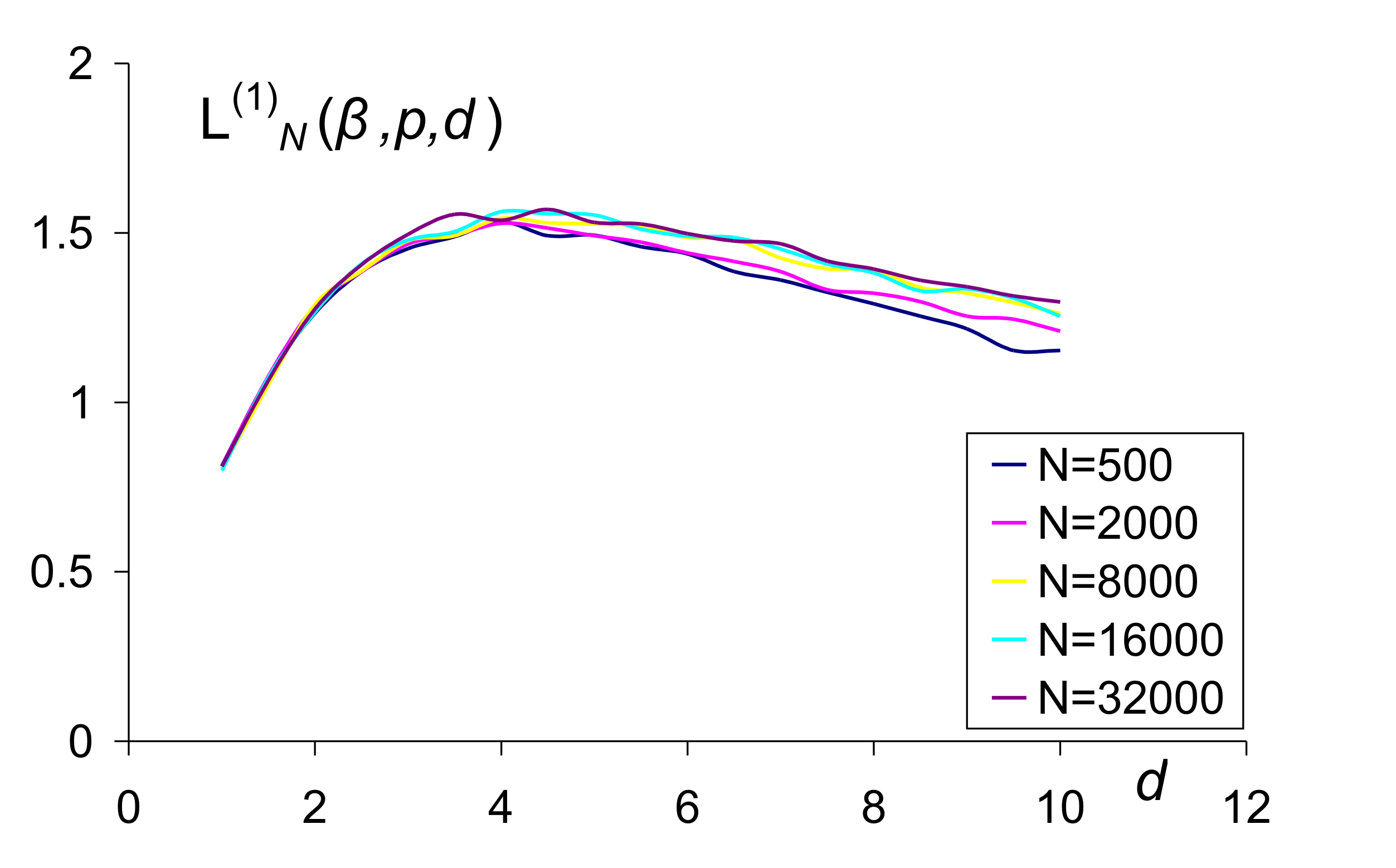}
\caption{Normalized regret for {\bf Algorithm 1}. $\beta=2.2$;
$p=0.5$; $N=500,2000,8000,16000,32000$. } \label{fig21}
\end{figure}
In this section, we provide a description of the MDA proposed
in~\cite{UN} for Bernoulli two-armed bandit problem. Note that the
idea of the mirror descent originates from~\cite{NU} for
multy-armed bandit.
\par
Let's introduce probability vectors
$\op_n=(p_{n}^{(1)},p_{n}^{(2)})$ s.t. $p_{n}^{(1)}\ge 0$,
$p_{n}^{(2)} \ge 0$, $p_{n}^{(1)}+p_{n}^{(2)}=1$, dual vectors
$\oz_n=(\zeta_{n}^{(1)},\zeta_{n}^{(2)})$ and stochastic gradient
vectors $\ou_n=(u_{n}^{(1)},u_{n}^{(2)})$. Gibbs distribution is
defined as follows
\begin{gather*}
    \oG_\beta(\oz)=\{S_\beta(\oz)\}^{-1}\left(e^{-\zeta^{(1)}/\beta},e^{-\zeta^{(2)}/\beta}    \right)
\end{gather*}
where
\begin{gather*}
    S_\beta(\oz)=e^{-\zeta^{(1)}/\beta}+e^{-\zeta^{(2)}/\beta}.
\end{gather*}

MDA for the two-armed bandit is defined recursively.


{\bf Algorithm 1.}

\vskip 3mm

\noindent \fbox{%
\parbox{8.5cm}{%

\begin{enumerate}
\item Fix some $\op_0$ and $\oz_0$.

\item For $n=1,2,\dots,N$:

\begin{enumerate}
\item Draw an action $y_n$ distributed as follows:

$\Pr\left( y_n=\ell\right)=p_{n-1}^{(\ell)}$, $\ell=1,2$;

\item Apply the action $y_n$ and get random income $\xi_n$
distributed as follows:
\begin{gather*}
\Pr\left( \xi_n=1|y_n=\ell\right)=p_\ell, \\
\Pr\left(\xi_n=0|y_n=\ell\right)=q_\ell,
\end{gather*}
$\ell=1,2$;

\item Compute the thresholded stochastic gradient
$\ou_{n}(\op_{n-1})$:

\parbox{5.5cm}{
\begin{gather*}
\ou_{n}(\op_{n-1})=\left\{
\begin{array}{l}
\left(\displaystyle{\frac{1-\xi_{n}}{p_{n-1}^{(1)}}},0\right), \ \mbox{if } \ y_{n}=1,\\
\left(0,\displaystyle{\frac{1-\xi_{n}}{p_{n-1}^{(2)}}}\right), \
\mbox{if } \ y_{n}=2;
\end{array}
\right.
\end{gather*}
}

\item Update the dual and probability vectors

$\oz_{n}=\oz_{n-1}+ \ou_n(\op_{n-1})$,

$\op_n=\oG_{\beta_n}(\oz_n)$;

\end{enumerate}

\end{enumerate}

}}

\vskip 3mm

\par
\begin{figure}[h]
\centering
\includegraphics[width=3.2in]{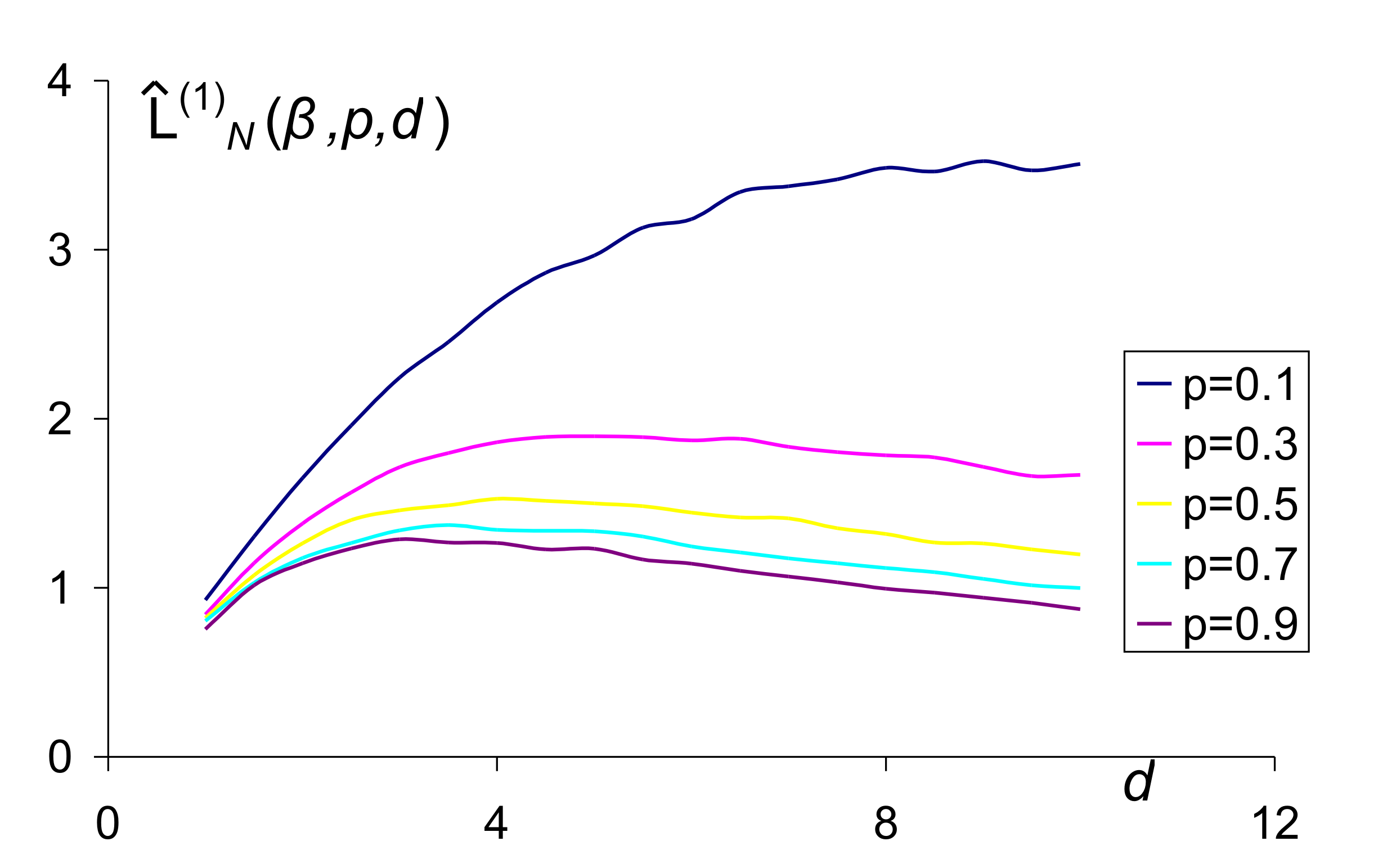}
\caption{Normalized regret for {\bf Algorithm 1}. $\beta=2.2$;
$N=2000$; $p=0.1,0.3,0.5,0.7,0.9$.} \label{fig22}
\end{figure}

Let's denote by $\Sigma_1$ the set of strategies described by the
MDA and by
\begin{equation}\label{a4}
R^{(1)}_N(\Theta)=\inf_{\Sigma_1} \sup_{\Theta} L_N(\sigma,\theta)
\end{equation}
corresponding minimax risk. The following theorem holds~\cite{UN}.
\par
\begin{theorem}\label{th3}
Consider {\bf Algorithm 1}. Let $\beta_n=\beta^*\times
\{D(n+1)\}^{1/2}$ with $\beta^*=(8/\log 2)^{1/2}\approx 3.397$,
$D=0.25$. Then for any horizon $N \ge 1$ for the minimax
risk~\eqref{a4} the estimate holds
\begin{gather}
\label{a5} R^{(1)}_N(\Theta) \le r^*_1 \{D (N+1)\}^{1/2}.
\end{gather}
with $r^*_1=4 (2 \log 2)^{1/2}\approx 4.710$.
\end{theorem}
\par
\begin{remark}
Our description of the algorithm differs from the original in some
details. The algorithm in~\cite{UN} is proposed for the problem of
minimization of the total expected income; it is done for
multi-armed bandit with arbitrary finite number of actions and
for 2nd moment rather than variance $D$ of incomes.
\end{remark}
\par
The estimate \eqref{a5} was obtained theoretically. It is
approximately $7.39$ times worse than the estimate \eqref{a3}.
However, it may be improved by Monte-Carlo simulations. To this
end, the following normalized regret was calculated:
\begin{gather*}
\hL_N^{(1)}(\beta,p,d)=(\hD N)^{-1/2}L_N(\sigma_N,\theta_N),
\end{gather*}
where $\theta_N$ and $d$ are related as $\theta=(p+d
(\hD/N)^{1/2},p-d (\hD/N)^{1/2})$, where $0<p<1$, $\hD=pq$,
$q=1-p$ and $\sigma_N$ stands for {\bf Algorithm 1} with
$\beta_n=\beta \{\hD(n+1)\}^{1/2}$. Here and below we put
$p_{0}^{(1)}=p_{0}^{(2)}=0.5$,
$\zeta_{0}^{(1)}=\zeta_{0}^{(2)}=0$. The number of Monte-Carlo
simulations is always 10000.
\par
On figure~\ref{fig21} we present $\hL_N^{(1)}(\beta,p,d)$
calculated for different horizons $N$ by Monte-Carlo simulations
if $\beta=2.2$, $p=0.5$ and $1 \le d \le 10$. Results are
presented for $N=500, 2000, 8000, 16000, 32000$. One can see that
$\hL_N^{(1)}(\beta,p,d)$ converges to some limiting function as $N
\to \infty$.
\par
One can guess that the limiting function $\hL_N^{(1)}(\beta,p,d)$
does not depend on $p$ if $0<p<1$ just like the results of
\cite{Koln1}, \cite{Koln2}. However, this is not the case for MDA.
On figure~\ref{fig22} we present $\hL_N^{(1)}(\beta,p,d)$
calculated by Monte-Carlo simulations if $\beta=2.2$, $N=2000$ and
$0 \le d \le 10$. Results are presented for $p=0.1, 0.3, 0.5, 0.7,
0.9$. One can see that the the lines are not similar and maximal
expected losses are attained for the smallest $p$.
\par
\begin{figure}[h]
\centering
\includegraphics[width=3.2in]{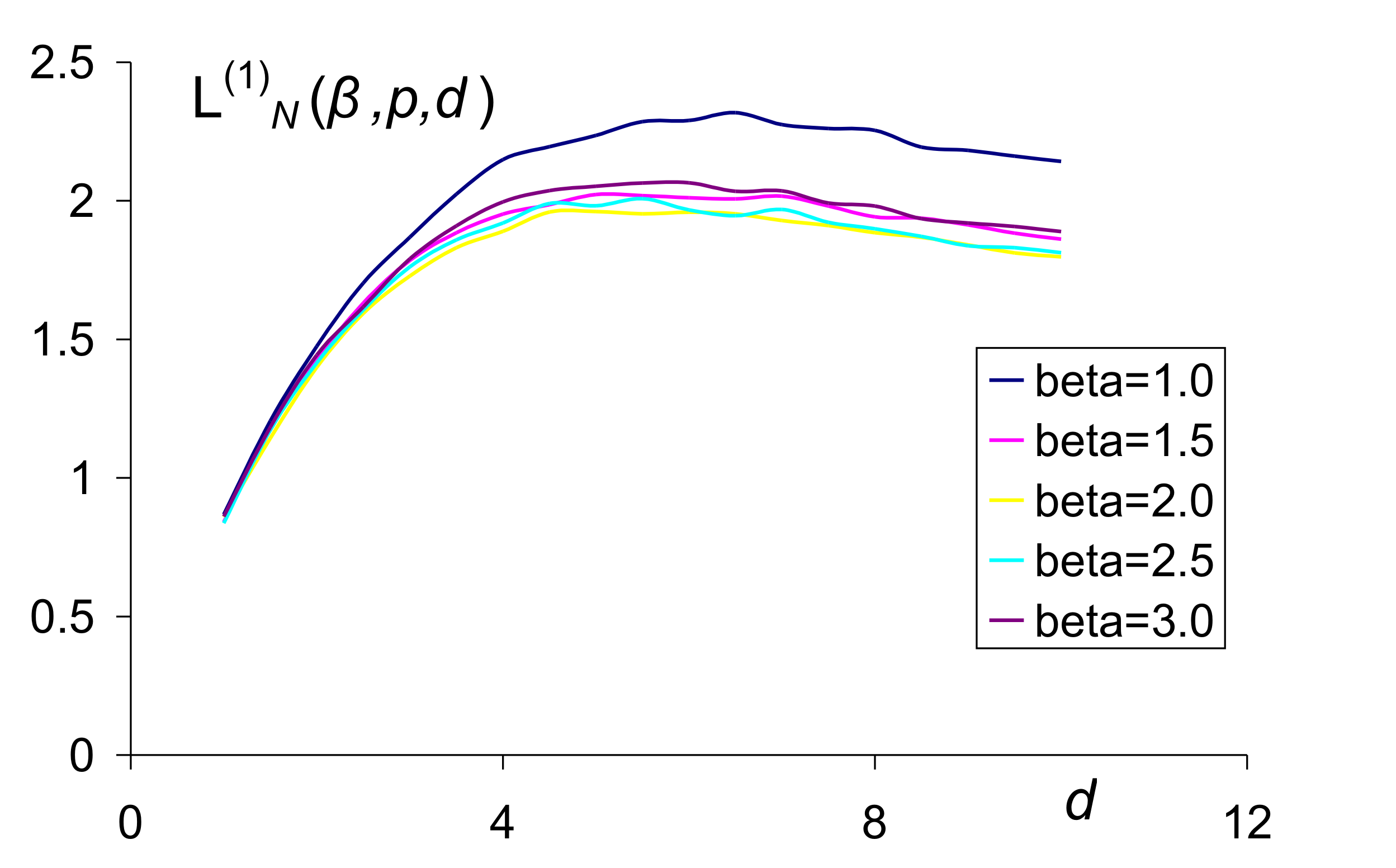}
\caption{Normalized regret for {\bf Algorithm 1}. $p=0.1$,
$N=2000$, $\beta=1.0, 1.5, 2.0, 2.5, 3.0$.} \label{fig23}
\end{figure}

Therefore, we calculate the following normalized regret
\begin{gather*}
\L_N^{(1)}(\beta,p,d)=(D N)^{-1/2}L_N(\sigma_N,\theta_N),
\end{gather*}
where $\theta_N$ and $d$ are related as $\theta=(p+d
(D/N)^{1/2},p-d (D/N)^{1/2})$ where $0<p<1$, $D=0.25$ and
$\sigma_N$ stands for {\bf Algorithm 1} with $\beta_n=\beta
\{D(n+1)\}^{1/2}$. First, we fix $p=0.1$ and calculate
$\L_N^{(1)}(\beta,p,d)$ by Monte-Carlo simulations if $N=2000$ and
$0 \le d \le 10$. Results are presented on figure~\ref{fig23} for
$\beta=1.0, 1.5, 2.0, 2.5, 2.0$. One can see that $\beta=2.0$ is
approximately optimal because it provides the least maximal
normalized regret in this case. More careful calculations give
that $\beta\approx 2.2$ is approximately optimal.
\par
\begin{figure}[h]
\centering
\includegraphics[width=3.2in]{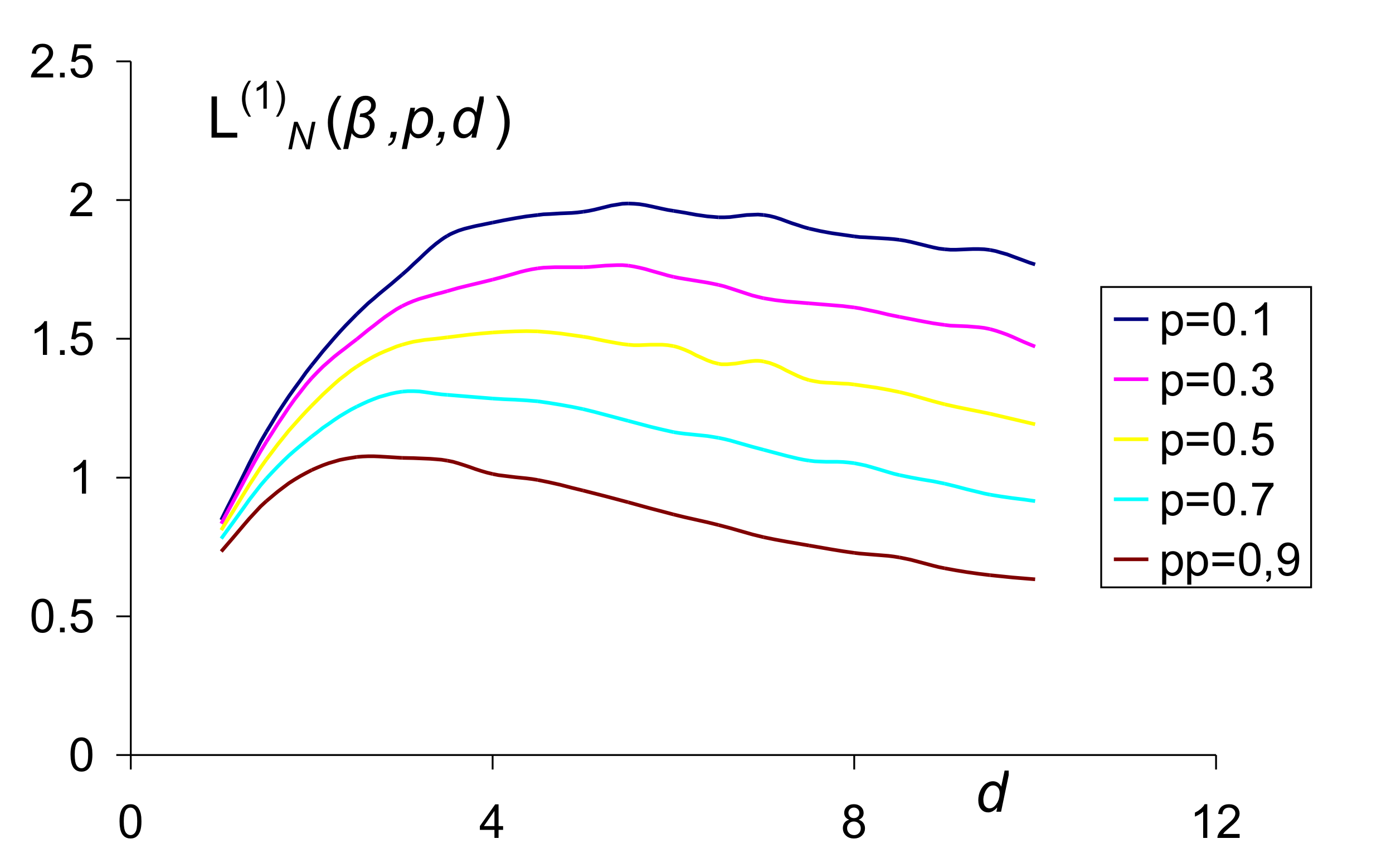}
\caption{Normalized regret for {\bf Algorithm 1} if $\beta=2.2$,
$N=2000$, $p=0.1, 0.3, 0.5, 0.7, 0.9$.} \label{fig24}
\end{figure}

Finally we calculate $\L_N^{(1)}(\beta,p,d)$ if $\beta\approx
2.2$, $N=2000$ for different $p$. Results are presented on
figure~\ref{fig24} for $p=0.1, 0.3, 0.5, 0.7, 0.9$. One can see
that maximal values of $\L_N^{(1)}(\beta,p,d)$ are attained if
$p=0.1$. Hence, the value $\beta \approx 2.2$ is approximately
optimal and
\begin{gather*}
r_1=\inf_{\beta>0}\max_{\footnotesize \begin{array}{l}
                         1 \le d\le 10,\\
                        0.1<p<0.9
                        \end{array}
                        }
\L_N^{(1)}(\beta,p,d) \approx 2.0.
\end{gather*}
This estimate is approximately 2.37 times better than theoretical
estimate \eqref{a5}.
\par

\section{Parallel Version of the MDA} \label{S3}

Parallel version of the MDA for Bernoulli two-armed bandit is also
defined recursively. Let's assume that $N=TM$ items of data are
given where $M$ describes the size of the packet of processed data
and $T$ is the number of processing stages. Let $\oM_t=(M_t^{(1)},
M_t^{(2)})$ be a vector s.t. $M_t^{(1)}>0$, $M_t^{(2)}>0$,
$M_t^{(1)}+M_t^{(2)}=M$. Denote by $[M_t^{(1)}]$, $[M_t^{(2)}]$
the closest integer numbers to $M_t^{(1)}$, $M_t^{(2)}$. We also
introduce a projection operator $\mP_\varrho(\op')=\op$ where
$\varrho>0$ and
\begin{gather*}
    \begin{array}{l}
    \op=\op', \ \mbox{if }  p'^{(1)}\ge \varrho, p'^{(2)}) \ge \varrho,\\
    p^{(1)}=\varrho, p^{(2)}=1-\varrho, \ \mbox{if } p'^{(1)}<
    \varrho,\\
    p^{(1)}=1-\varrho, p^{(2)}=\varrho, \ \mbox{if } p'^{(2)}<
    \varrho.
    \end{array}
\end{gather*}
The following parallel version of the MDA assigns $[M_t^{(1)}]$,
$[M_t^{(2)}]$ in $t$-th packet in proportion to $p_{t-1}^{(1)}$,
$p_{t-1}^{(2)}$ and then applies the first and the second actions
$[M_t^{(1)}]$ and $[M_t^{(2)}]$ times respectively.

\begin{figure}[h]
\centering
\includegraphics[width=3.2in]{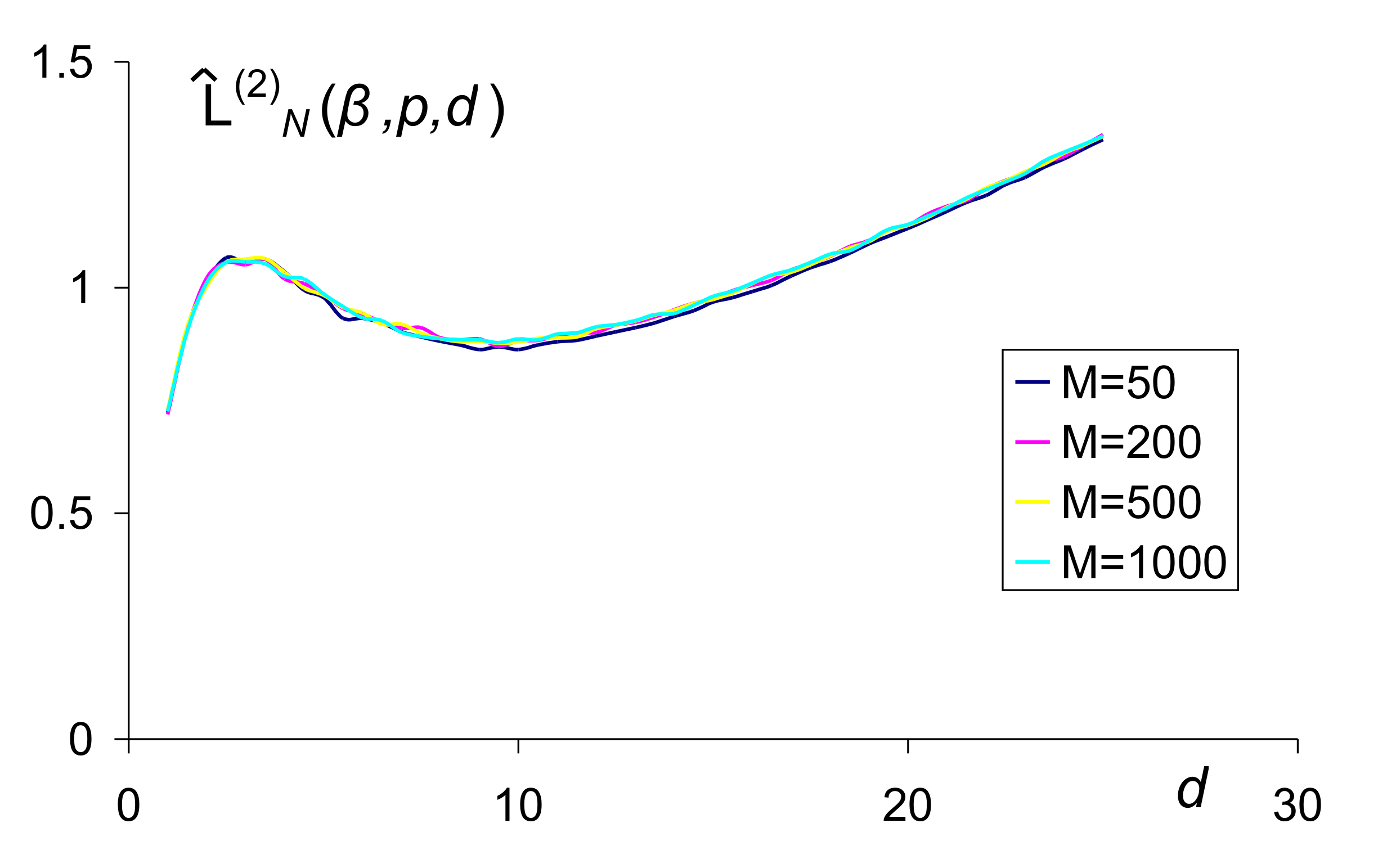}
\caption{Normalized regret for {\bf Algorithm~2} if $\beta=1.0$;
$\varrho=0.02$; $T=100$; $M=50, 200, 500, 1000$.} \label{fig31}
\end{figure}

\vskip 1mm

{\bf Algorithm 2.}

\vskip 3mm

\noindent \fbox{%
\parbox{8.5cm}{%

\begin{enumerate}
\item Fix some $\op_0$ and $\oz_0$.

\item For $t=1,2,\dots,T$:

\begin{enumerate}
\item Let $M_t^{(\ell)}=p_{t-1}^{(\ell)}\times M$, $\ell=1,2$;

\item For $n \in [(t-1)M+1,tM]$ apply the $\ell$-th action
$[M_t^{(\ell)}]$ times getting
random income \\
$\eta_t^{(\ell)}=\sum (1-\xi_n|y_n=\ell)$,\\
where
\begin{gather*}
\Pr\left( \xi_n=1|y_n=\ell\right)=p_\ell, \\
\Pr\left(\xi_n=0|y_n=\ell\right)=q_\ell,
\end{gather*}
$\ell=1,2$;

\item Compute the thresholded stochastic gradient
$\ou_{t}(\op_{t-1})$:

\parbox{5.5cm}{
\begin{gather*}
\ou_{t}(\op_{t-1})=
\left(\displaystyle{\frac{\eta_{t}^{(1)}}{p_{t-1}^{(1)}}},
\displaystyle{\frac{\eta_{t}^{(2)}}{p_{t-1}^{(2)}}}\right),
\end{gather*}
}

\item Update the dual and probability vectors

$\oz_{t}=\oz_{t-1}+ \ou_t(\op_{t-1})$,

$\op'_t=G_{\beta_{t}}(\oz_t)$,

$\op_t=\mP_\varrho (\op'_t)$;

\end{enumerate}

\end{enumerate}

}}
\vskip 3mm
\par

The projection operator is used because $M_t^{(\ell)}<1$  if
corresponding $p_{t-1}^{(\ell)}$ is small enough. We define the
normalized regret
\begin{gather*}
\hL_N^{(2)}(\beta,p,d)=(\hD N)^{-1/2}L_N(\sigma_N,\theta_N)
\end{gather*}
where $\theta_N=(p+d(\hD/N)^{1/2},p-d(\hD/N)^{1/2})$, and
$\sigma_N$ stands for {\bf Algorithm 2} with $\beta_t=\beta (\hD
M(t+0.5))^{1/2}$, $0<p<1$, $q=1-p$, $\hD=pq$.
\par
Note that if $M$ is large enough then according to the central
limit theorem all $\{\eta^{(\ell)}_t\}$ are approximately normally
distributed with parameters
\begin{gather*}
\mE (\eta^{(\ell)}_t)=q_\ell [M_t^{(\ell)}],\
\V(\eta^{(\ell)}_t)=\hD [M_t^{(\ell)}],
\end{gather*}
$\ell=1,2$. For normally distributed incomes, it is convenient to
set a control problem with a continuous time. Namely, if the
$\ell$-th action is applied for a duration $\Delta M$, which is
not obligatory integer, then it generates normally distributed
income with mathematical expectation $q_\ell \times \Delta M$ and
variance $\hD \times \Delta M$ and independent from all previously
obtained incomes. Let's present corresponding modified algorithm.

\vskip 3mm

{\bf Algorithm 3.}

\vskip 3mm

\noindent \fbox{%
\parbox{8.5cm}{%

\begin{enumerate}
\item Fix some $\op_0$ and $\oz_0$.

\item For $t=1,2,\dots,T$:

\begin{enumerate}
\item Let $M_t^{(\ell)}=p_{t-1}^{(\ell)}\times M$, $\ell=1,2$;

\item Apply the $\ell$-th action for a duration $M_t^{(\ell)}$
getting
random income $\eta_t^{(\ell)}$ which is normally distributed with \\
$\mathbb{E}(\eta_t^{(\ell)})=q_\ell \times M_t^{(\ell)}$,\\
$\V(\eta_t^{(\ell)})=\hD \times M_t^{(\ell)}$, $\ell=1,2$ and
independent from all previous incomes;

\item Compute the thresholded stochastic gradient
$\ou_{t}(\op_{t-1})$:

\parbox{5.5cm}{
\begin{gather*}
\ou_{t}(\op_{t-1})=
\left(\displaystyle{\frac{\eta_{t}^{(1)}}{p_{t-1}^{(1)}}},
\displaystyle{\frac{\eta_{t}^{(2)}}{p_{t-1}^{(2)}}}\right),
\end{gather*}
}

\item Update the dual and probability vectors

$\oz_{t}=\oz_{t-1}+ \ou_t(\op_{t-1})$,

$\op'_t=G_{\beta_{t}}(\oz_t)$,

$\op_t=\mP_\varrho (\op'_t)$;

\end{enumerate}

\end{enumerate}

}}

\vskip 3mm

We define the normalized regret
\begin{gather*}
\hL_N^{(3)}(\beta,p,d)=(\hD N)^{-1/2}L_N(\sigma_N,\theta_N)
\end{gather*}
where $\theta_N=(p+d(\hD/N)^{1/2},p-d(\hD/N)^{1/2})$, and
$\sigma_N$ stands for {\bf Algorithm 3} with $\beta_t=\beta (\hD
M(t+0.5))^{1/2}$, $0<p<1$, $q=1-p$, $\hD=pq$.
\par
\begin{theorem}\label{th4}
Consider  {\bf Algorithm 3} with a fixed number $T$ of packet
processing stages and arbitrary $\hD>0$. Then normalized loss
function $\mathrm{L}_N^{(3)}(\beta,p,d)$ does not depend on $N$,
$M$, $p$, $\hD$.
\end{theorem}
\par
{\bf Proof.} Consider $q_\ell=w_\ell (\hD/N)^{1/2}$, $\ell=1,2$,
so as $q_1-q_2=-2d (\hD/N)^{1/2}$, $d>0$. Let
$X^{(\ell)}_t(D^{(\ell)}_t)$ denote independent normally
distributed random variables s.t.
\begin{gather*}
    \mE(X^{(\ell)}_t(D^{(\ell)}_t))=0,\
   \V (X^{(\ell)}_t(D^{(\ell)}_t))= D^{(\ell)}_t,
\end{gather*}
$\ell=1,2$; $t=1,2,\dots,T$. Then
\begin{gather*}
\eta^{(\ell)}_t=w_\ell (\hD/N)^{1/2} M_{t}^{(\ell)} +
X^{(\ell)}_t(\hD M_{t}^{(\ell)})\\ =\eps w_\ell (\hD N)^{1/2}
p_{t-1}^{(\ell)} + X^{(\ell)}_t( \hD N \eps p_{t-1}^{(\ell)}).
\end{gather*}
Here $\eps=M/N$. Next, we obtain
\begin{gather*}
\frac{\eta^{(\ell)}_t}{p_{t-1}^{(\ell)}} =\eps w_\ell (\hD
N)^{1/2} + X^{(\ell)}_t\left(\frac{\hD  N \eps}
{p_{t-1}^{(\ell)}}\right).
\end{gather*}
So,
\begin{gather} \label{a6}
\zeta^{(\ell)}_{t}=\tau w_\ell (\hD N)^{1/2}+ \sum_{i=1}^t
X^{(\ell)}_i\left( \frac{\hD N \eps } {p_{i-1}^{(\ell)}}\right),
\end{gather}
where $\tau=t/T$. Recall that $\beta_t=\beta (\hD
M(t+0.5))^{1/2}$. Hence
\begin{gather*}
\frac{\zeta^{(\ell)}_{t}}{\beta_{t}}=\frac{1}{\beta}\left(\frac{\tau
w_\ell}{(\tau+ 0.5\eps)^{1/2}} + Y^{(\ell)}(t)\right)
\end{gather*}
with
\begin{gather*}
Y^{(\ell)}(t)=\sum_{i=1}^t X^{(\ell)}_i\left( \frac{\eps }
{(\tau+0.5 \eps) p_{i-1}^{(\ell)}}\right).
\end{gather*}
So,
\begin{gather}\label{a7}
\begin{array}{c}
\displaystyle{\frac{\zeta^{(1)}_{t}-\zeta^{(2)}_{t}}{\beta_{t}}}
=-\frac{ 2d \tau} {\beta (\tau+ 0.5\eps)^{1/2}}\\
+ \displaystyle{\frac{1}{\beta}} \left(
Y^{(1)}(t)-Y^{(2)}(t)\right).
\end{array}
\end{gather}
Note that \eqref{a7} depends only on parameters of packet
processing. Since
\begin{gather} \label{a8}
p_{t-1}'^{(1)}=
\frac{\exp\left(\displaystyle{-\frac{\zeta^{(1)}_{t}-\zeta^{(2)}_{t}}{\beta_{t}}}\right)}
{\exp\left(\displaystyle{-\frac{\zeta^{(1)}_{t}-\zeta^{(2)}_{t}}{\beta_{t}}}\right)+1}
\end{gather}
and $p_{t-1}'^{(2)}=1-p_{t-1}'^{(1)}$ then all
$\{p_{t-1}^{(\ell)}\}$ depend only on the parameters of packet
processing. The regret $\mathrm{L}_N^{(3)}(\beta,p,d)$ can be
expressed as follows
\begin{gather}\label{a9}
\begin{array}{c}
\mathrm{L}_N^{(3)}(\beta,p,d)\\ =(\hD N)^{-1/2}(p_1-p_2)
\displaystyle{\sum_{t=1}^T} M \
\mathbb{E} \left(p_{t-1}^{(2)}\right)\\
=2d \displaystyle{\sum_{t=1}^T} \eps \
\mathbb{E}\left(p_{t-1}^{(2)}\right).
\end{array}
\end{gather}
This expression does not depend on $N$, $M$, $p$, $\hD$ but only
on the parameters of packet processing. This proves
theorem~\ref{th4}.
\par
\begin{figure}[h]
\centering
\includegraphics[width=3.2in]{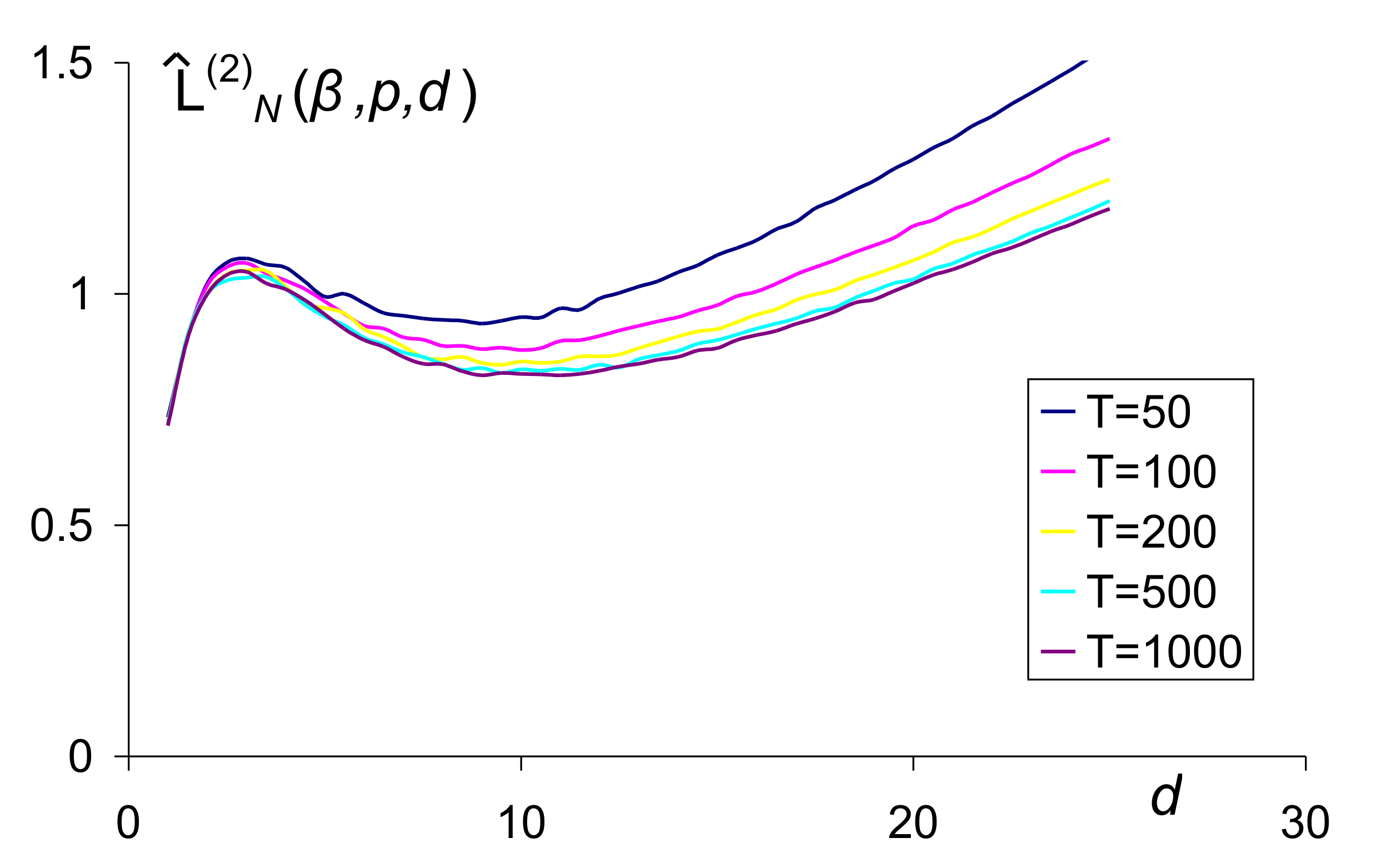}
\caption{Normalized regret for {\bf Algorithm~2} if $\beta=1.0$;
$\varrho=0.02$; $M=100$; $T=50, 100, 200, 500, 1000$.}
\label{fig32}
\end{figure}

\begin{figure}[h]
\centering
\includegraphics[width=3.2in]{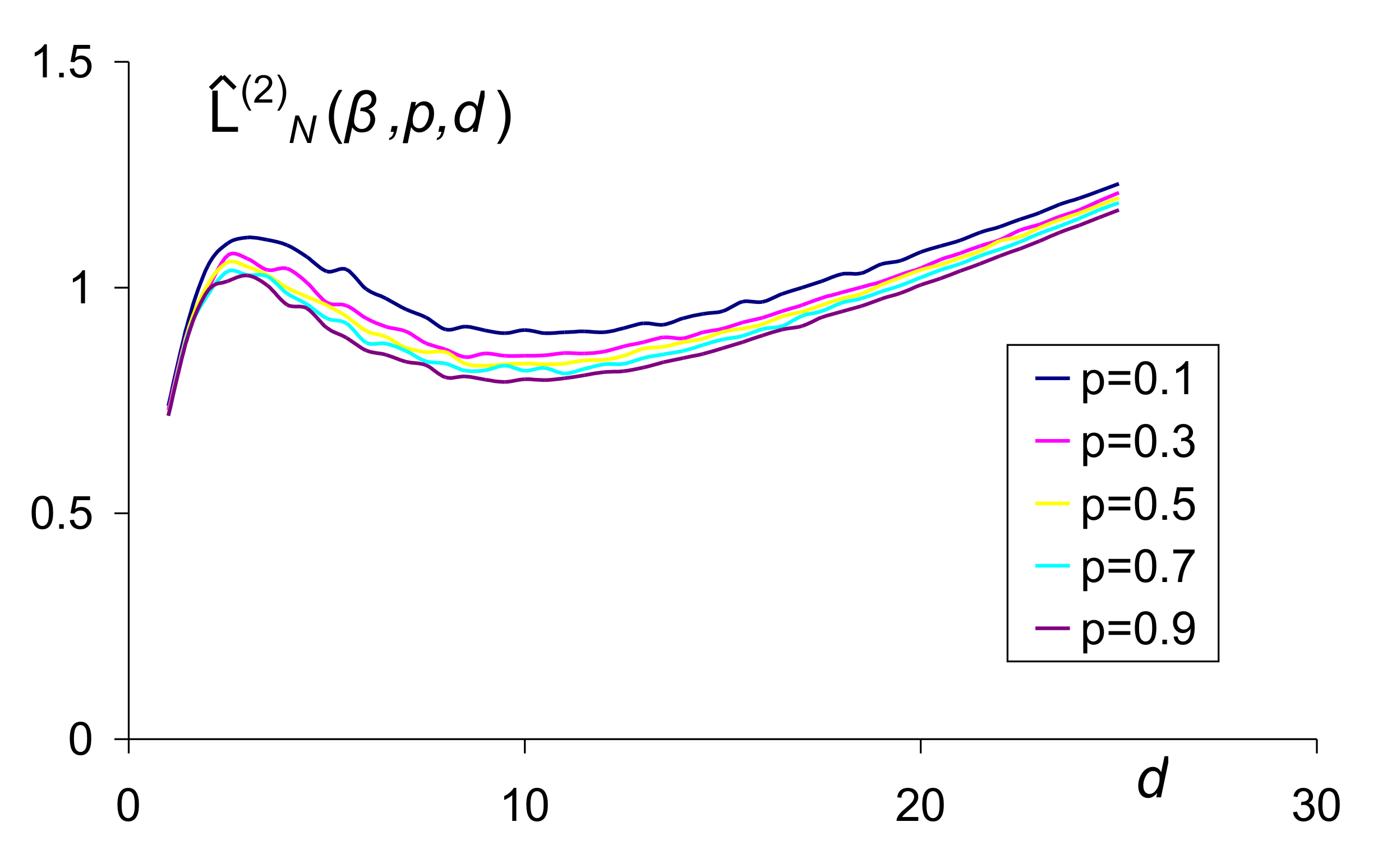}
\caption{Normalized regret for {\bf Algorithm~2} if $\beta=1.0$;
$\varrho=0.02$; $M=100$; $T=500$; $p=0.1, 0.3, 0.5, 0.7, 0.9$.}
\label{fig33}
\end{figure}

\begin{figure}[h]
\centering
\includegraphics[width=3.2in]{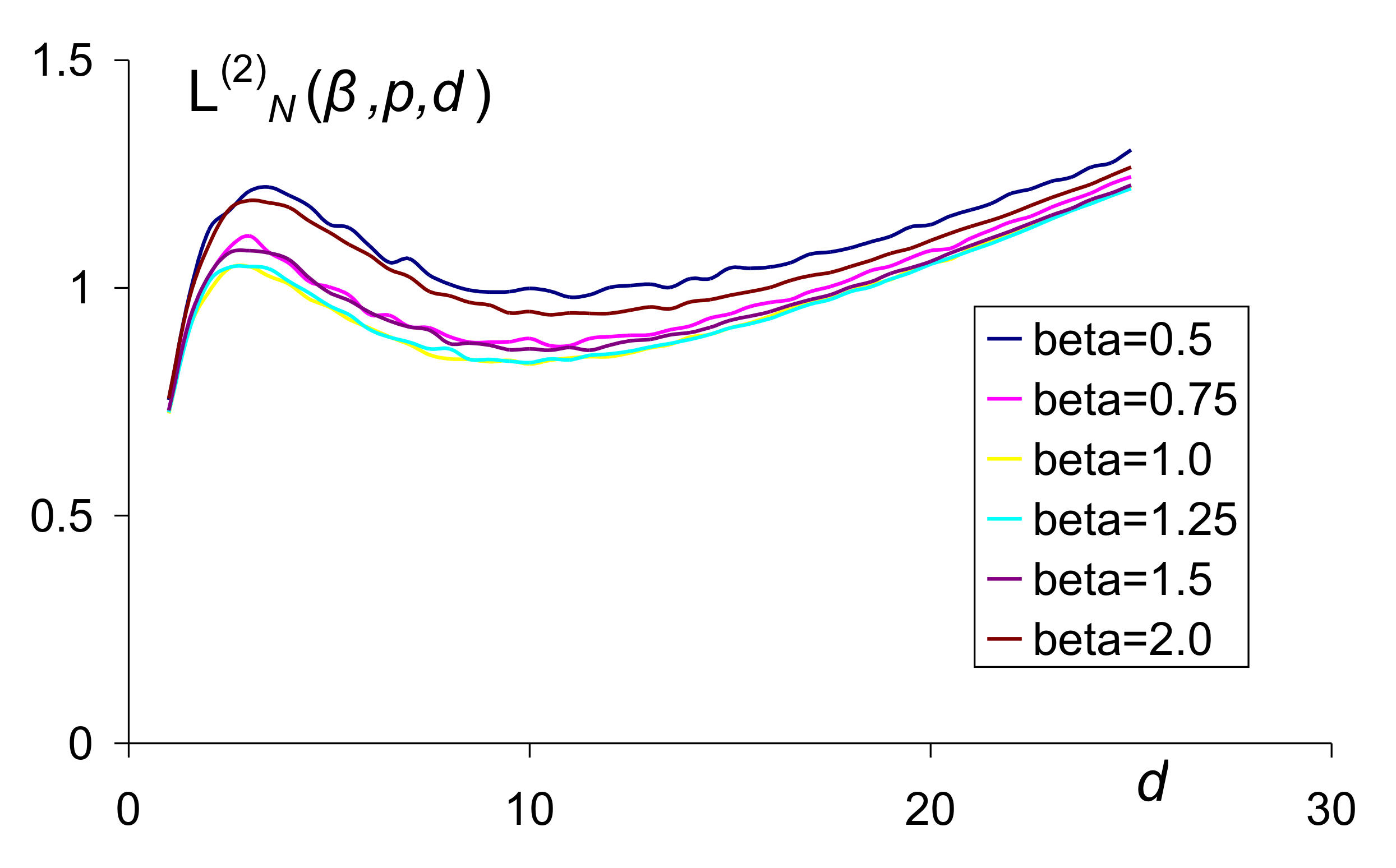}
\caption{Normalized regret for {\bf Algorithm~2} if
$\varrho=0.02$; $M=100$; $T=300$; $p=0.5$; $\beta=0.5, 0.75, 1.0,
1.25, 1.5, 2.0$.} \label{fig34}
\end{figure}

\begin{figure}[h]
\centering
\includegraphics[width=3.2in]{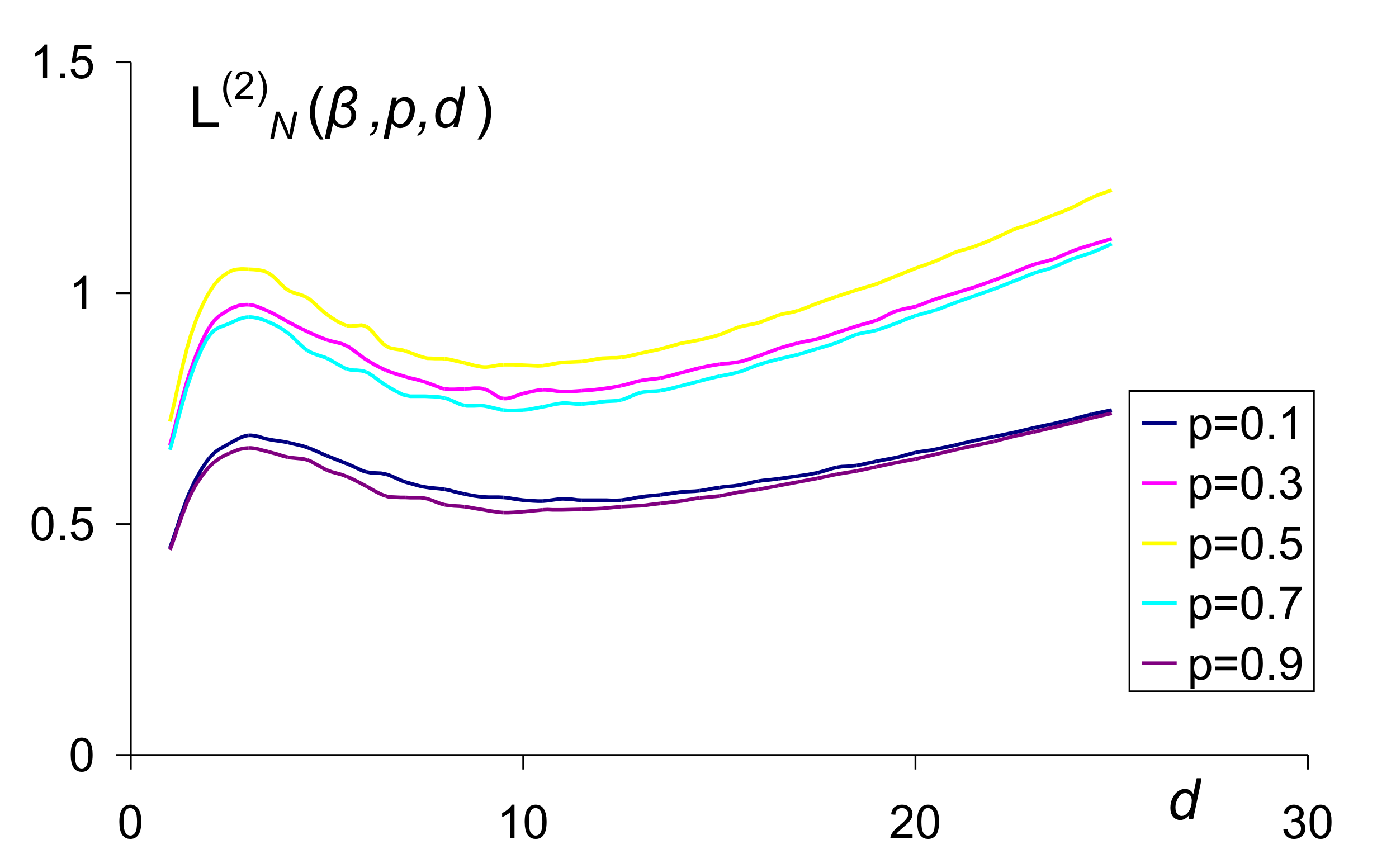}
\caption{Normalized regret for {\bf Algorithm~2} if $\beta=1.0$;
$\varrho=0.02$; $M=100$; $T=300$; $p=0.1, 0.3, 0.5, 0.7, 0.9$.}
\label{fig35}
\end{figure}
\par
\begin{remark}\label{r5}
Denote by $\{\heta^{(\ell)}_t\}$ incomes corresponding to {\bf
Algorithm 2}. Let $\d^{(\ell)}_t$ be independent random variables
s.t.
\begin{gather*}
    \mE(\d^{(\ell)}_t)=O(M^{-1}),\
   \V (\d^{(\ell)}_t)= O(M^{-1}).
\end{gather*}
Then $\heta^{(\ell)}_t=\eta^{(\ell)}_t(1+\d^{(\ell)}_t)$,
$\ell=1,2$; $t=1,2,\dots,T$, because deviations
$\{\d^{(\ell)}_t\}$ are caused by nonintegral values
$M_{t}^{(\ell)}$ in {\bf Algorithm 2}. Since the number of stages
$T$ is finite and fixed it means that normalized loss function
$\mathrm{L}_N^{(2)}(\beta,p,d)$ is close to
$\mathrm{L}_N^{(3)}(\beta,p,d)$ in \eqref{a9} if $M$ is large
enough.
\end{remark}
\par
\begin{remark}
It seems very likely that {\bf Algorithm 3} converges as $N \to
\infty$ and $\eps=M/N \to 0$. Let's put
$\varsigma_N^{(\ell)}(\tau)=(\hD N)^{-1/2} \zeta^{(\ell)}_{t}$,
$P_N^{(\ell)}(\tau)=p^{(\ell)}_{t-1}$, $\ell=1,2$ where
$\tau=t/T$. Let $W_{\ell}(\tau)$, $\ell=1,2$ be independent Wiener
processes. Denote by $\varsigma^{(\ell)}(\tau)$,
$P^{(\ell)}(\tau)$ corresponding limiting random processes and
$\mathrm{L}^{(3)}(\beta,p,d)$ a limiting normalized regret as $N
\to \infty$ and $\eps=M/N \to 0$. Using \eqref{a6}, \eqref{a8} and
\eqref{a9} a limiting description may be presented as
\begin{gather*}
d \varsigma^{(\ell)}(\tau)= w_\ell d\tau+
(P^{(\ell)}(\tau))^{-1/2} d W_\ell(\tau),
\end{gather*}
\begin{gather*}
P'^{(1)}(\tau)=
\frac{\exp\left(\displaystyle{-\frac{\varsigma^{(1)}(\tau)-\varsigma^{(2)}(\tau)}
{\beta \tau^{1/2}}}\right)}
{\exp\left(\displaystyle{-\frac{\varsigma^{(1)}(\tau)-\varsigma^{(2)}(\tau)}{\beta
\tau^{1/2}}}\right)+1},
\end{gather*}
$P'^{(2)}(\tau)=1-P'^{(1)}(\tau)$,  $\oP(\tau)=\mP\{\oP'(\tau)\}$,
$\tau \in [0,1]$. Initial conditions are given by
\begin{gather*}
\varsigma^{(1)}(0)=\varsigma^{(2)}(0)=0.
\end{gather*}
Then limiting normalized regret is equal to
\begin{gather*}
\mathrm{L}^{(3)}(\beta,p,d)= 2d \int_{0}^1 \mathbb{E}
\left(P^{(2)}(\tau)\right) d \tau .
\end{gather*}
However, we do not have a rigorous proof of this result.
\end{remark}
\par
In view of remark~\ref{r5} we present simulations of {\bf
Algorithm 2} but expect to observe the properties of {\bf
Algorithm 3}. On figure~\ref{fig31} we present
$\hL_N^{(2)}(\beta,p,d)$ calculated for different sizes of packet
$M$ by Monte-Carlo simulations if $\beta=1.0$, $p=0.5$,
$\varrho=0.02$, $T=100$ and $1 \le d \le 25$. Results are
presented for $M=50, 200, 500, 1000$. One can see that
$\hL_N^{(2)}(\beta,p,d)$ is almost independent from the size of
packet $M$.
\par
On figure~\ref{fig32} we present $\hL_N^{(2)}(\beta,p,d)$
calculated for different $T$ by Monte-Carlo simulations if
$\beta=1.0$, $p=0.5$, $\varrho=0.02$, $M=100$ and $1 \le d \le
25$. Results are presented for $T=50, 100, 200, 500, 1000$. One
can see that $\hL_N^{(2)}(\beta,p,d)$ converges as $T \to \infty$.
\par
According to theorem \ref{th4} the limiting function
$\hL_N^{(2)}(\beta,p,d)$ as $N \to \infty$ does not depend on $p$
if $0<p<1$. On figure~\ref{fig33} we present
$\hL_N^{(2)}(\beta,p,d)$ calculated by Monte-Carlo simulations if
$\beta=1.0$, $M=100$, $T=500$, $\varrho=0.02$ and $1 \le d \le
25$. Results are presented for $p=0.1, 0.3, 0.5, 0.7, 0.9$. One
can see that these lines are close to each other.
\par
Then we determine the optimal $\beta$. We define the normalized
regret
\begin{gather*}
\L_N^{(2)}(\beta,p,d)=(D N)^{-1/2}L_N(\sigma_N,\theta_N)
\end{gather*}
where $\theta_N=(p+d(D/N)^{1/2},p-d(D/N)^{1/2})$ and $\sigma_N$
stands for {\bf Algorithm 2} with $\beta_t=\beta (D
M(t+0.5))^{1/2}$, $0<p<1$, $D=0.25$. First, we fix $p=0.5$ and
calculate $\L_N^{(2)}(\beta,p,d)$ by Monte-Carlo simulations if
$M=100$, $T=300$, $\varrho=0.02$ and $0 \le d \le 25$. Results are
presented on figure~\ref{fig34} for $\beta=0.5, 0.75, 1.25, 1.5,
2.0$. One can see that $\beta=1.0$ is approximately optimal
because it provides the least maximal normalized regret
$\L_N^{(2)}(\beta,p,d) < 1.05$ if $d<20$.
\par
Finally we calculate $\L_N^{(2)}(\beta,p,d)$ if $\beta= 1.0$,
$M=100$, $T=300$, $\varrho=0.02$ and $0 \le d \le 25$. Results are
presented on figure~\ref{fig35} for $p=0.1, 0.3, 0.5, 0.7, 0.9$.
One can see that maximal values of $\L_N^{(2)}(\beta,p,d)$ are
attained if $p=0.5$. Hence, the value $\beta = 1.0$ is
approximately optimal and
\begin{gather*}
r_2=\inf_{\beta>0}\max_{\footnotesize \begin{array}{l}
                         1 \le d\le 20,\\
                        0.1<p<0.9
                        \end{array}
                        }
\L_N^{(2)}(\beta,p,d) \approx 1.1.
\end{gather*}
This estimate is even better than $r_1$. However, it is attained
for close $p_1$, $p_2$ because $\varrho>0$.

\section{Another Parallel Version of the MDA} \label{S4}
\par
Let's consider now another parallel version of the MDA which
behaves closely to the ordinary version.

\vskip 3mm

{\bf Algorithm 4.}

\vskip 3mm

\noindent \fbox{%
\parbox{8.5cm}{%

\begin{enumerate}
\item Fix some $\op_0$ and $\oz_0$.

\item For $t=1,2,\dots,T$:

\begin{enumerate}
\item
 \begin{enumerate}
 \item Put $\chi^{(1)}_{t}=\chi^{(2)}_{t}=0$.

 \item For $n=(t-1)\times M+1,\dots,t\times M$:

  \begin{enumerate}
   \item Draw an action $y_n$ distributed as follows:\\
   $\Pr\left( y_n=\ell\right)=p_{t-1}^{(\ell)}$, $\ell=1,2$;

   \item Apply the action $y_n$, get random income $\xi_n$
   distributed as follows:
   \begin{gather*}
   \Pr\left( \xi_n=1|y_n=\ell\right)=p_\ell, \\
   \Pr\left(\xi_n=0|y_n=\ell\right)=q_\ell,
   \end{gather*}
   and update:
   \begin{gather*}
   \chi^{(\ell)}_{t}\leftarrow \chi^{(\ell)}_{t}+(1-\xi_n)
   \end{gather*}
   if $y_n=\ell$,\ $\ell=1,2$;.
  \end{enumerate}
 \end{enumerate}

\item Compute the thresholded stochastic gradient
$\ou_{t}(\op_{t-1})$:

\parbox{5.5cm}{
\begin{gather*}
\ou_{t}(\op_{t-1})=
\left(\displaystyle{\frac{\chi^{(1)}_{t}}{p_{t-1}^{(1)}}},
\frac{\chi^{(2)}_{t}}{p_{t-1}^{(2)}}  \right),
\end{gather*}
}

\item Update the dual and probability vectors

$\oz_{t}=\oz_{t-1}+ \ou_t(\op_{t-1})$,

$\op_t=\oG_{\beta_t}(\oz_t)$;

\end{enumerate}

\end{enumerate}

}}

\vskip 3mm

\par
\begin{figure}[h]
\centering
\includegraphics[width=3.2in]{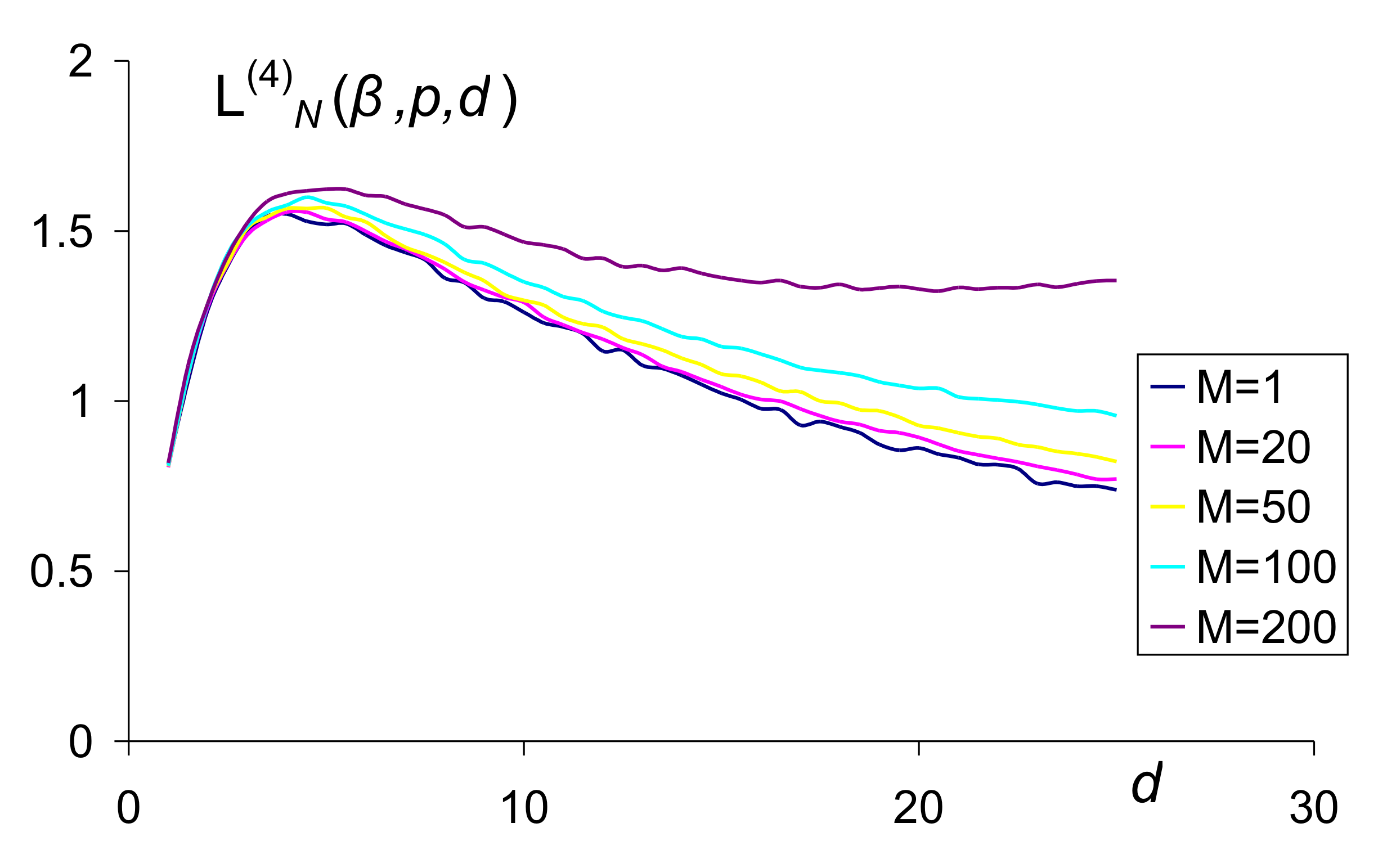}
\caption{Normalized regret for {\bf Algorithm~4} if $\beta=2.2$;
$N=10000$; $p=0.5$; $M=100$; $M=1, 20, 50, 100, 200$.}
\label{fig41}
\end{figure}
\par
It is straightforward to check that
\begin{gather*}
\chi^{(\ell)}_{t}=\sum_{n=(t-1)M+1}^{tM} (1-\xi^{(\ell)}_n),
\end{gather*}
where $\{\xi^{(\ell)}_n\}$ are i.i.d. variables distributed as
\begin{gather*}
\Pr\{(1-\xi^{(\ell)}_n)=1\}= p_{t-1}^{(\ell)} q_\ell, \\
\Pr\{(1-\xi^{(\ell)}_n)=0\}= 1-p_{t-1}^{(\ell)} q_\ell
\end{gather*}
and $\{\xi^{(1)}_n\}$, $\{\xi^{(2)}_n\}$ are independent from each
other. Hence
\begin{gather}
\begin{array}{c} \label{a10}
\mE \chi^{(\ell)}_{t}=Mp_{t-1}^{(\ell)} q_\ell,\\
\V  \chi^{(\ell)}_{t}=Mp_{t-1}^{(\ell)} q_\ell (1-p_{t-1}^{(\ell)}
q_\ell), \ \ell=1,2.
\end{array}
\end{gather}
\par
Let's put $p_1=p+d(D/N)^{1/2}$, $p_2=p-d(D/N)^{1/2}$, $D=0.25$ and
$q_\ell=1-p_\ell$, $\ell=1,2$. Note that if $N$ is large enough
then distributions of $\chi^{(\ell)}_{t}$, $\ell=1,2$ are close to
gaussian ones. We define the normalized regret
\begin{gather*}
\L_N^{(4)}(\beta,p,d)=(D N)^{-1/2}L_N(\sigma_N,\theta_N),
\end{gather*}
where $\theta_N=(p+d(D/N)^{1/2},p-d(D/N)^{1/2})$, and $\sigma_N$
stands for {\bf Algorithm 4} with $\beta_t=\beta
(DM(t+0.5))^{1/2}$.
\par
\begin{theorem}\label{th5}
Consider  {\bf Algorithm 4} with a fixed number $T$ of packet
processing stages. Assume that $\chi^{(\ell)}_t$, $\ell=1,2$ are
normally distributed with mathematical expectations and variances
given by~\eqref{a10}. Then asymptotically (as $N \to \infty$)
normalized loss function $\L_N^{(4)}(\beta,p,d)$ does not depend
on $N$, $M$, but does depend on on $q$, $D$ and parameters of
packet processing.
\end{theorem}
\par
{\bf Proof.} Let's put $q_\ell=w_\ell (D/N)^{1/2}$, $\ell=1,2$, so
that $q_1-q_2=-2d (D/N)^{1/2}$, $d>0$. Let
$X^{(\ell)}_t(D^{(\ell)}_t)$ denote independent normally
distributed random variables s.t.
\begin{gather*}
    \mathbb{E}(X^{(\ell)}_t(D^{(\ell)}_t))=0,\
   \V (X^{(\ell)}_t(D^{(\ell)}_t))= D^{(\ell)}_t,
\end{gather*}
$\ell=1,2$; $t=1,2,\dots,T$. Then
\begin{gather*}
\chi^{(\ell)}_t=w_\ell (D/N)^{1/2} M p_{t-1}^{(\ell)} +
X^{(\ell)}_t(Mp_{t-1}^{(\ell)} q_\ell (1-p_{t-1}^{(\ell)} q_\ell))\\
=\eps w_\ell (D N)^{1/2} p_{t-1}^{(\ell)} + X^{(\ell)}_t( \eps
Np_{t-1}^{(\ell)} q_\ell (1-p_{t-1}^{(\ell)} q_\ell)).
\end{gather*}
Here $\eps=M/N$. Next, we obtain
\begin{gather*}
\frac{\chi^{(\ell)}_t}{p_{t-1}^{(\ell)}} =\eps w_\ell (D N)^{1/2}
+ X^{(\ell)}_t\left(\frac{\eps N  q_\ell (1-p_{t-1}^{(\ell)}
q_\ell)} {p_{t-1}^{(\ell)}}\right).
\end{gather*}
So,
\begin{gather} \label{a11}
\zeta^{(\ell)}_{t}=\tau w_\ell (D N)^{1/2}+ \sum_{i=1}^t
X^{(\ell)}_i\left(\frac{\eps N  q_\ell (1-p_{t-1}^{(\ell)}
q_\ell)} {p_{t-1}^{(\ell)}}\right),
\end{gather}
where $\tau=t/T$. Recall that $\beta_t=\beta
(DM(t+0.5))^{1/2}=\beta (DN(\tau+0.5\eps))^{1/2}$. Hence
\begin{gather*}
\frac{\zeta^{(\ell)}_{t}}{\beta_{t}}=\frac{1}{\beta}\left(\frac{\tau
w_\ell}{(\tau+ 0.5\eps)^{1/2}} + Y^{(\ell)}_N(t)\right)
\end{gather*}
with
\begin{gather*}
Y^{(\ell)}_N(t)=\sum_{i=1}^t X^{(\ell)}_i\left( \frac{\eps q_\ell
(1-p_{t-1}^{(\ell)} q_\ell)} {D(\tau+0.5 \eps)
p_{i-1}^{(\ell)}}\right).
\end{gather*}
Note that $q_\ell \to q$ as $M \to \infty$, $\ell=1,2$. Hence,
$Y^{(\ell)}_N(t) \to Y^{(\ell)}(t)$ as $M \to \infty$ where
\begin{gather*}
Y^{(\ell)}(t)= \sum_{i=1}^t X^{(\ell)}_i\left( \frac{\eps q
(1-p_{t-1}^{(\ell)} q)} {D(\tau+0.5 \eps)
p_{i-1}^{(\ell)}}\right).
\end{gather*}
As the number of stages $T$ is fixed, then asymptotically (as $M
\to \infty$)
\begin{gather}\label{a12}
\begin{array}{c}
\displaystyle{\frac{\zeta^{(1)}_{t}-\zeta^{(2)}_{t}}{\beta_{t}}}
=-\frac{ 2d \tau} {\beta (\tau+ 0.5\eps)^{1/2}}\\
+ \displaystyle{\frac{1}{\beta}} \left(
Y^{(1)}(t)-Y^{(2)}(t)\right).
\end{array}
\end{gather}
Note that \eqref{a12} depends only on $q$ , $D$ and $T$. Since
\begin{gather} \label{a13}
p_{t-1}^{(1)}=
\frac{\exp\left(\displaystyle{-\frac{\zeta^{(1)}_{t}-\zeta^{(2)}_{t}}{\beta_{t}}}\right)}
{\exp\left(\displaystyle{-\frac{\zeta^{(1)}_{t}-\zeta^{(2)}_{t}}{\beta_{t}}}\right)+1}
\end{gather}
and $p_{t-1}^{(2)}=1-p_{t-1}^{(1)}$ then all
$\{p_{t-1}^{(\ell)}\}$ depend only on $q$, $D$ and $T$. The regret
$\L_N^{(4)}(\beta,p,d)$ can be expressed as follows
\begin{gather}\label{a14}
\begin{array}{c}
\mathrm{L}_N^{(4)}(\beta,p,d)\\ =(D N)^{-1/2}(p_1-p_2)
\displaystyle{\sum_{t=1}^T} M \
\mathbb{E} \left(p_{t-1}^{(2)}\right)\\
=2d \displaystyle{\sum_{t=1}^T} \eps \
\mathbb{E}\left(p_{t-1}^{(2)}\right).
\end{array}
\end{gather}
This expression does not depend on $N$, $M$,  but does depend on
$q$, $D$ and $T$. This proves theorem~\ref{th4}.
\par
\begin{remark}
Like in section~\ref{S3}, we present a limiting description of
{\bf Algorithm 4} as $M \to \infty$ and $\eps=M/N \to 0$. Let's
put $\varsigma_N^{(\ell)}(\tau)=(D N)^{-1/2} \zeta^{(\ell)}_{t}$,
$P_N^{(\ell)}(\tau)=p^{(\ell)}_{t-1}$, $\ell=1,2$ where
$\tau=t/T$. Let $W_{\ell}(\tau)$, $\ell=1,2$ be independent Wiener
processes. Denote by $\varsigma^{(\ell)}(\tau)$,
$P^{(\ell)}(\tau)$ corresponding limiting random processes as $N
\to \infty$ and $\eps=M/N \to 0$ and $\L^{(4)}(\beta,p,d)$ the
limit of $\L^{(4)}_N(\beta,p,d)$. Using \eqref{a11}, \eqref{a13}
and \eqref{a14} a limiting description may be presented as
\begin{gather*}
d \varsigma^{(\ell)}(\tau)= w_\ell d\tau+ \left(\frac{q(1-q
P^{(\ell)}(\tau))}{DP^{(\ell)}(\tau)}\right)^{1/2} d W_\ell(\tau),
\end{gather*}
\begin{gather*}
P^{(1)}(\tau)=
\frac{\exp\left(\displaystyle{-\frac{\varsigma^{(1)}(\tau)-\varsigma^{(2)}(\tau)}
{\beta \tau^{1/2}}}\right)}
{\exp\left(\displaystyle{-\frac{\varsigma^{(1)}(\tau)-\varsigma^{(2)}(\tau)}{\beta
\tau^{1/2}}}\right)+1},
\end{gather*}
where $w_1-w_2=-2d$, $P^{(2)}(\tau)=1-P^{(1)}(\tau)$, $\tau \in
[0,1]$. Initial conditions are given by
\begin{gather*}
\varsigma^{(1)}(0)=\varsigma^{(2)}(0)=0.
\end{gather*}
Then limiting normalized regret is equal to
\begin{gather*}
\mathrm{L}^{(4)}(\beta,p,d)= 2d \int_{0}^1 \mathbb{E}
\left(P^{(2)}(\tau)\right) d \tau .
\end{gather*}
We do not have a rigorous proof of this result as well.
\end{remark}
\par
On figure~\ref{fig41} we present $\L_N^{(4)}(\beta,p,d)$
calculated for different $M$ by Monte-Carlo simulations if
$\beta=2.2$, $p=0.5$, $M=10000$ and $1 \le d \le 25$. Results are
presented for $M=20, 50, 100, 200$ (accordingly $T=500, 200, 100,
50$). The case $M=1$ corresponds to ordinary MDA and
$\L_N^{(1)}(\beta,p,d)$. One can see that $\L_N^{(4)}(\beta,p,d)$
is close to $\L_N^{(1)}(\beta,p,d)$ if $T$ is large enough.
\par
\section{Combined Algorithms} \label{S5}
\par
One can see on figure~\ref{fig32} and figure~\ref{fig41} that
larger sizes of packets correspond to larger sizes of normalized
regret if $d$ is large enough. It is caused by equal applications
of both actions to initial packet. To avoid this effect of initial
packet processing one can take initial packets of smaller sizes.
The simplest decision is to use the ordinary algorithm at initial
short stage and then to switch to parallel algorithm.
\par

\begin{figure}[h]
\centering
\includegraphics[width=3.2in]{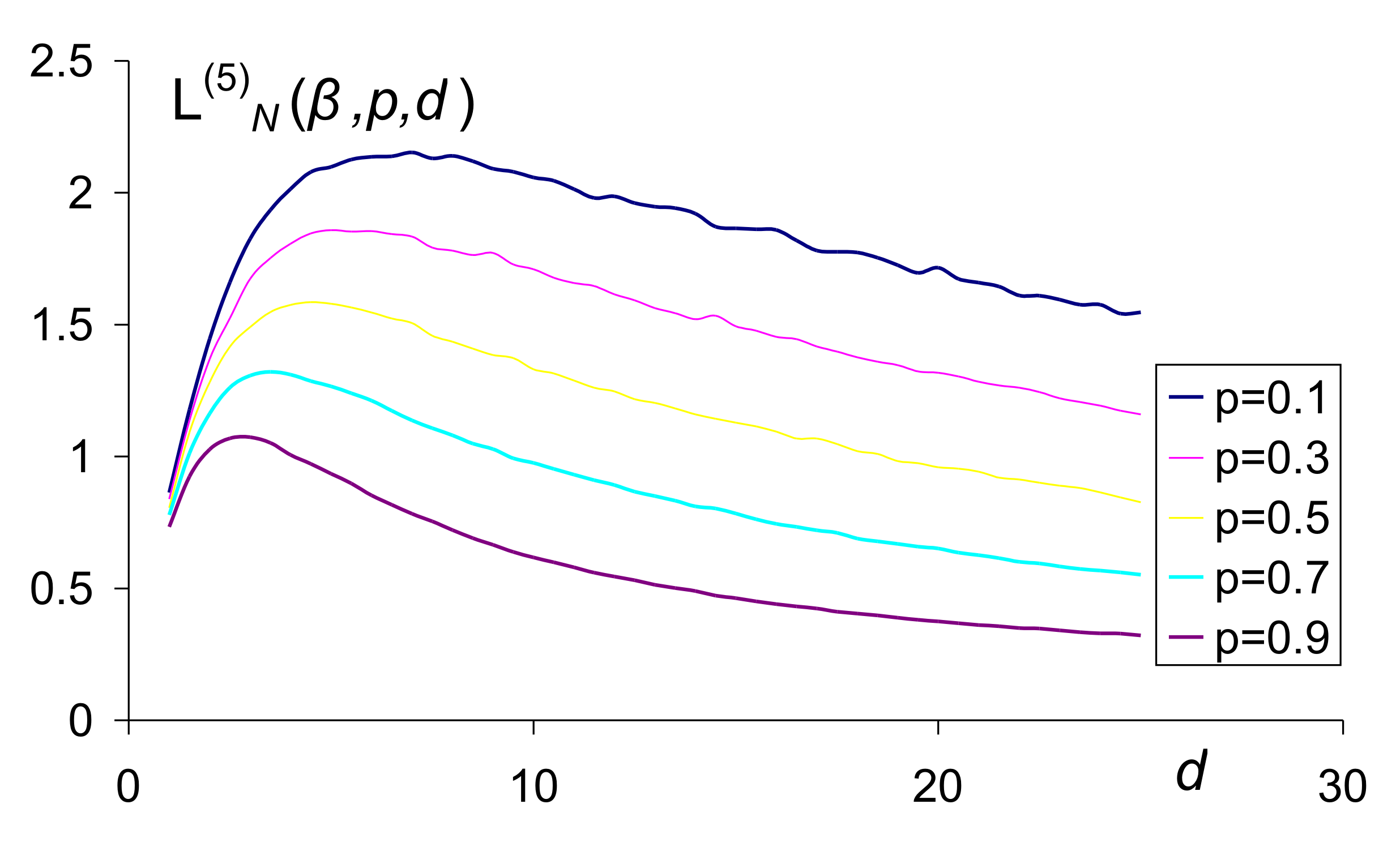}
\caption{Normalized regret for {\bf Algorithm~5} if $\beta=2.2$;
$N=20000$; $M_0=600$; $M=200$; $p=0.1, 0.3, 0.5, 0.7, 0.9$.}
\label{fig51}
\end{figure}
\par

\begin{figure}[h]
\centering
\includegraphics[width=3.2in]{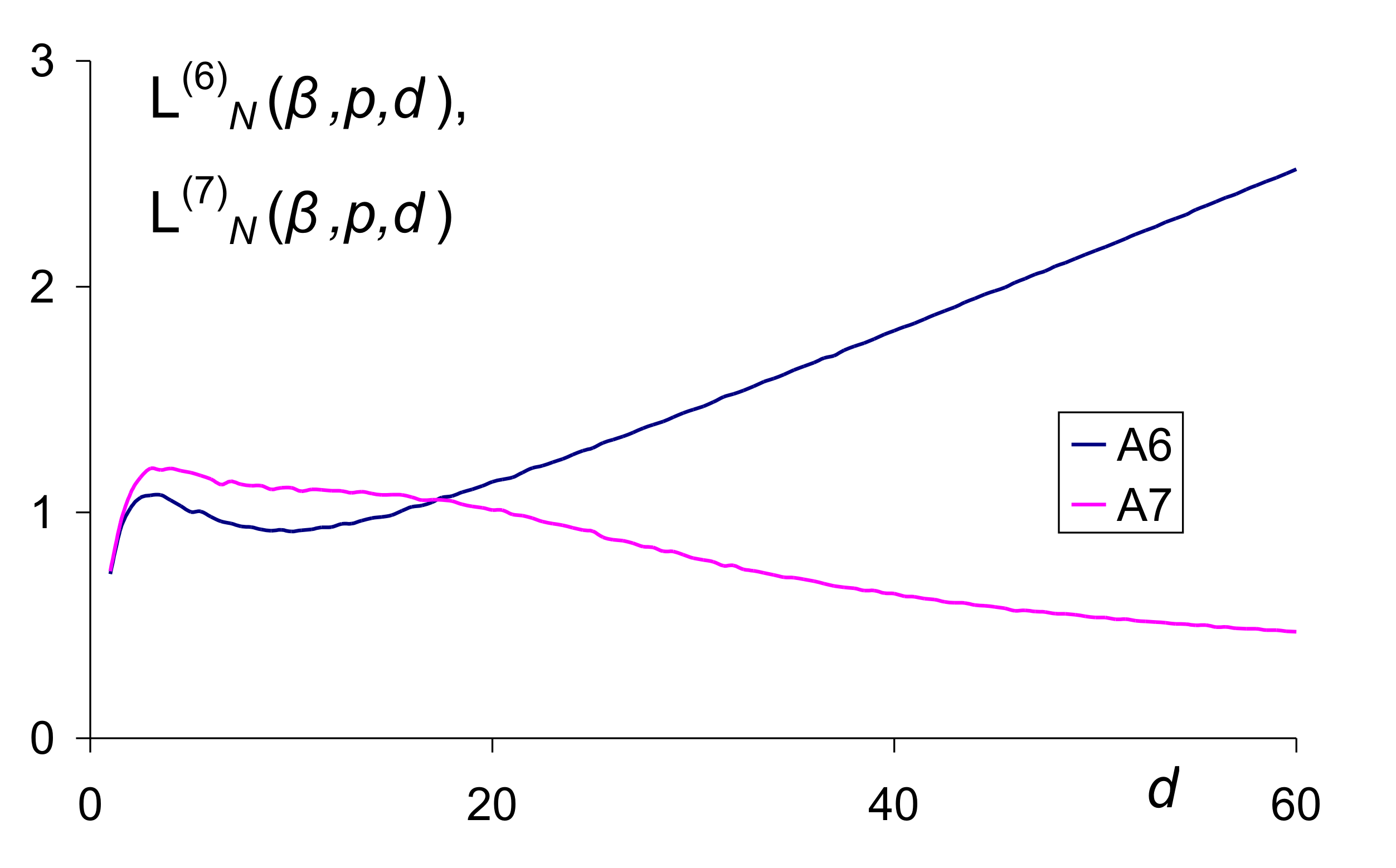}
\caption{Comparative results for {\bf Algorithm~6} and {\bf
Algorithm~7} if $\beta_1=2.2$; $\beta_2=1.0$; $N=30000$; $p=0.5$;
$M_0=900$; $M=300$, $\varrho=0.02$, $\kappa=0.2$.} \label{fig52}
\end{figure}
\par

\begin{figure}[h]
\centering
\includegraphics[width=3.2in]{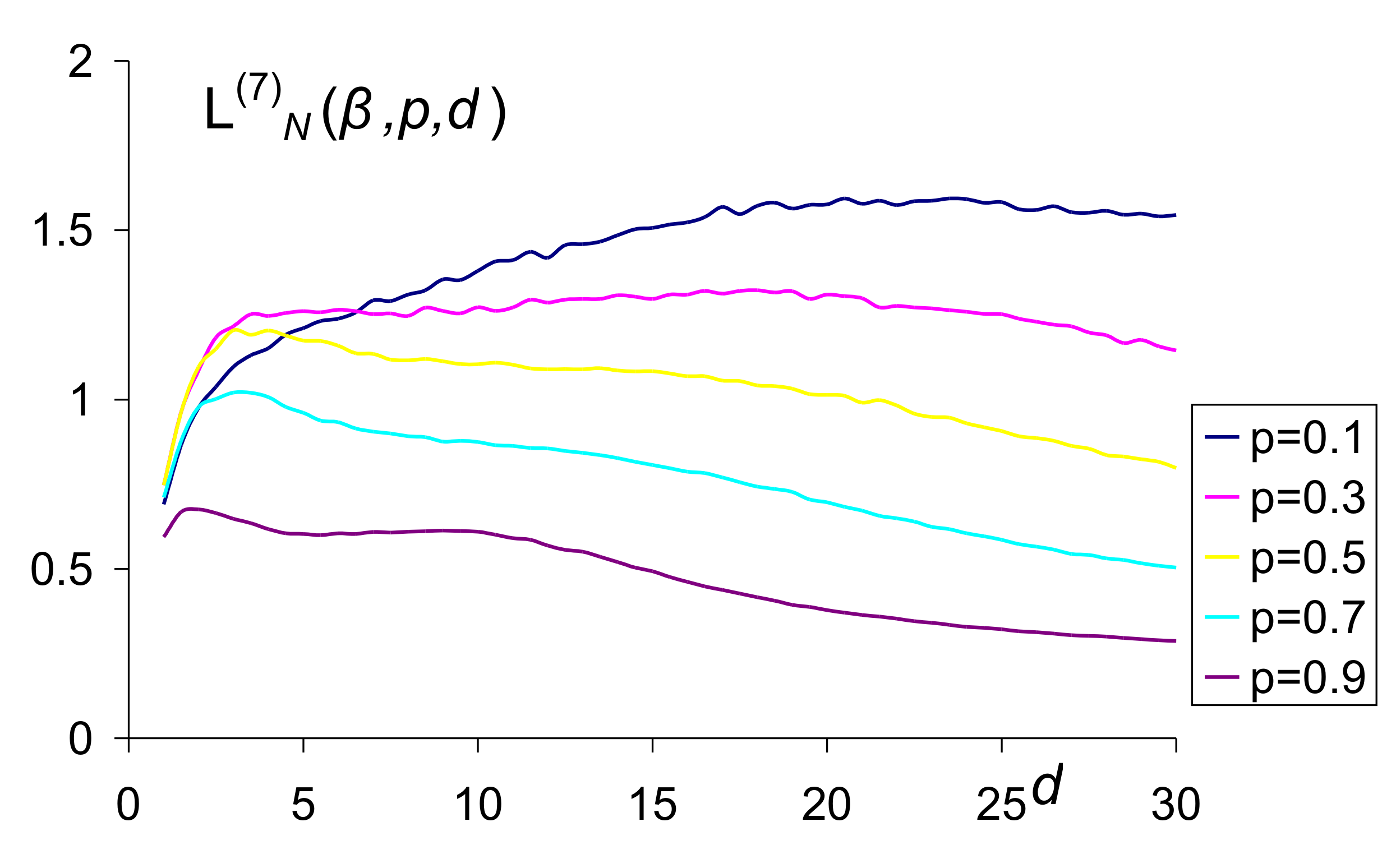}
\caption{Normalized regret for {\bf Algorithm~7} if $\beta_1=2.2$;
$\beta_2=1.0$; $N=30000$; $M_0=900$; $M=300$, $\varrho=0.02$,
$\kappa=0.2$; $p=0.1, 0.3, 0.5, 0.7, 0.9$.} \label{fig53}
\end{figure}

First, we combine {\bf Algorithm 1} and {\bf Algorithm 4} as
follows.
\par
{\bf Algorithm 5.}

\vskip 3mm
\noindent \fbox{%
\parbox{8.5cm}{%

\begin{enumerate}
\item Apply {\bf Algorithm 1} at the initial horizon \\
$n=1,\dots,M_0$. Get $\op_{M_0}$ and $\oz_{M_0}$.

\item Apply {\bf Algorithm 4} at the residual horizon
$n=M_0+1,\dots,N$.

\end{enumerate}

}} \vskip 3mm

On figure~\ref{fig51} we present $\L_N^{(5)}(\beta,p,d)$
calculated by Monte-Carlo simulations  for {\bf Algorithm~5} if
$\beta=2.2$; $N=20000$; $M_0=600$; $M=200$; $p=0.1, 0.3, 0.5, 0.7,
0.9$. One can see that results are close to those presented on
figure~\ref{fig24} for ordinary MDA.
\par
To take advantage of {\bf Algorithm 2}, we combine {\bf Algorithm
1} and {\bf Algorithm 2} as follows.
\par
{\bf Algorithm 6.}

\vskip 3mm
\noindent \fbox{%
\parbox{8.5cm}{%

\begin{enumerate}
\item Apply {\bf Algorithm 1} at the initial horizon \\
$n=1,\dots,M_0$ with $\beta=\beta_1$. Get $\op_{M_0}$ and
$\oz_{M_0}$.

\item Apply {\bf Algorithm 2} at the residual horizon
$n=M_0+1,\dots,N$ with $\beta=\beta_2$.

\end{enumerate}

}} \vskip 3mm

However, {\bf Algorithm 2} provides large normalized regret if $d$
is large enough (see figure~\ref{fig31}, \ref{fig32}) because it
applies both actions with probabilities no less than $\varrho$.
Therefore we consider the following combined algorithm.

\par
{\bf Algorithm 7.}

\vskip 3mm
\noindent \fbox{%
\parbox{8.5cm}{%

\begin{enumerate}
\item Apply {\bf Algorithm 1} at the initial horizon \\
$n=1,\dots,M_0$ with $\beta=\beta_1$. Get $\op_{M_0}$ and
$\oz_{M_0}$.

\item If $\min (p_{M_0}^{(1)},p_{M_0}^{(2)})<\kappa$\\
then apply {\bf Algorithm 4} with $\beta=\beta_1$\\
else apply {\bf Algorithm 2} with $\beta=\beta_2$\\
at the residual horizon $n=M_0+1,\dots,N$.

\end{enumerate}

}} \vskip 3mm

If $\kappa$ is appropriately chosen this algorithm for small $d$
switches mostly to {\bf Algorithm 2}. For large $d$ it switches
mostly to {\bf Algorithm 4}. On figure~\ref{fig52} we present
comparative results for $\L_N^{(6)}(\beta,p,d)$ and
$\L_N^{(7)}(\beta,p,d)$ if $\beta_1=2.2$; $\beta_2=1.0$;
$N=30000$; $p=0.5$; $M_0=900$; $M=300$, $\varrho=0.02$,
$\kappa=0.2$. One can see that $\L_N^{(7)}(\beta,p,d)$ does not
grow for large $d$.
\par
Finally, on figure~\ref{fig53} we present normalized regret for
{\bf Algorithm~7} if $\beta_1=2.2$; $\beta_2=1.0$; $N=30000$;
$M_0=900$; $M=300$, $\varrho=0.02$, $\kappa=0.2$; $p=0.1, 0.3,
0.5, 0.7, 0.9$. Maximal values of $\L_N^{(7)}(\beta,p,d)$ are
larger than those for $\L_N^{(2)}(\beta,p,d)$ but smaller than
those for $\L_N^{(5)}(\beta,p,d)$.
\par

\section{Conclusion} \label{Con}

Two parallel versions of the mirror descent algorithm (MDA) for
the two-armed bandit problem are proposed. The usage of parallel
versions of the MDA ensures that total time of data processing
depends mostly on the number of packets but not on the total
number of data. Monte-Carlo simulations show that maximal expected
losses for parallel versions are not more than for ordinary
version which processes data one-by-one. However, it is true only
for close mathematical expectations. For distant mathematical
expectations the effect of the first packet processing, when both
actions are equally applied, causes significant losses if the size
of the packet is large enough. This effect may be avoided if at
the sufficiently short initial stage the ordinary mirror descent
algorithm is used and then switched to the parallel version.



%

\end{document}